\documentclass[11pt,onecolumn]{article}
\usepackage{algorithm, algorithmic, setspace}
\usepackage{graphicx} 

\usepackage{amsmath} 
\usepackage{amssymb}  
\usepackage{amsthm}
\usepackage{amsfonts}
\usepackage{dsfont}
\usepackage{color}
\usepackage{subfig}
\usepackage{booktabs}
\usepackage{subfig}
\usepackage{pgfplots}
\usepackage{tikz}
\usetikzlibrary{arrows}
\usetikzlibrary{intersections}
\usepackage[english]{babel}
\usepackage[latin1]{inputenc}

\newcommand{\fun}{\mathcal{F}}
\newcommand{\lag}{\mathcal{L}}
\newcommand{\gun}{\mathcal{H}}

\newcommand{\soft}{\mathbb{S}}

\newcommand{\X}{\mathcal{X}}

\newcommand{\argmin}[1]{\underset{#1}{\mathrm{argmin\,}}}

\newcommand{\xmin}{x^{\star}}
\newcommand{\xtrue}{\widetilde{x}}
\newcommand{\bmin}{\beta^{\star}}
\newcommand{\btrue}{\widetilde{\beta}}
\newcommand{\alf}{\mathcal{A}}

\newcommand{\R}{\mathds{R}}
\newcommand{\N}{\mathds{N}}
\newcommand{\I}{\mathds{1}} 
\newcommand{\Z}{\mathds{Z}}

\newtheorem{lemma}{Lemma}
\newtheorem{proposition}{Proposition}
\newtheorem{theorem}{Theorem}
\newtheorem{corollary}{Corollary}
\newtheorem{remark}{Remark}
\newtheorem{definition}{Definition}
\newtheorem{assumption}{Assumption}

\begin{document}
\title{Non-convex Lasso-kind approach to compressed sensing for finite-valued signals}
\author{Sophie M. Fosson}
\maketitle
\begin{abstract}In this paper, we bring together two trends that have recently emerged in sparse signal recovery: the problem of sparse signals that stem from finite alphabets and the techniques that introduce concave penalties. Specifically, we show that using a minimax concave penalty (MCP) the recovery of finite-valued sparse signals is enhanced with respect to Lasso, in terms of estimation accuracy, number of necessary measurements, and run time. {\textcolor{black}{We focus on problems where sparse signals can be recovered from few linear measurements, as stated in compressed sensing theory.}} We start by proposing a Lasso-kind functional with MCP, whose minimum is the desired signal in the noise-free case, under null space conditions. We analyze its robustness to noise as well. We then propose an efficient ADMM-based algorithm to search the minimum. The algorithm is proven to converge to the set of stationary points, and its performance is evaluated through numerical experiments, both on randomly generated data and on a real localization problem. Furthermore, in the noise-free case, it is possible to check the exactness of the solution, and we test a version of the algorithm that exploits this fact to look for the right signal.
\end{abstract}

\section{Introduction}
In the last decade, the development of compressed sensing (CS, \cite{don06, fou11}) has brought a novel interest on sparse signals, that is, signals that have few non-zero components or that can be represented by few non-zero components in certain bases. Sparse signals are ubiquitous in diverse  applications. CS has established a new paradigm  by stating that sparse signals can be recovered from few non-adaptive linear measurements.  In this setting, a $k$-sparse signal $x\in\R^n$ can be recovered from the compressed measurements $y=Ax$ (possibly corrupted by noise) with $A\in\R^{m,n}$, $m<n$, under the assumptions on $A$ analyzed in the CS theory.  
\subsection{Finite-valued sparse signals}
In the extensive CS literature, some recent work is dedicated to the subcase of \emph{finite-valued} signals, \emph{i.e.}, $x\in\alf^n$ where $\alf$ is a {\textcolor{black}{known alphabet, that is, a finite set of symbols}}. This is a problem encountered in a number of sparse/CS applications, such as digital image recovery \cite{bio14}, security \cite{bio14sec}, digital communications \cite{spa15, ili12}, and discrete control signal design \cite{bem99}. In many localization problems \cite{fen09, bay15}, the localization area is split into cells and the goal is to verify which cells are occupied or not, the number of occupied cells being generally much smaller than the total: this can be interpreted as the recovery of a binary sparse signal in $\{0,1\}^n$. A binary framework is also present in spectrum sensing in cognitive radio networks and wideband spectrum sensing \cite{baz10,zen11,axe12, rom13}, where the goal is to detect if users are active or not in a given spectrum band. Furthermore, in many other applications, ranging from opinion polls to sensors data, intrinsically discrete or quantized data are envisaged \cite{das13}.

Even though the sub-problem of finite-valued signals can be approached by classical CS, it makes more sense to try to exploit the prior knowledge of $\alf$  to improve the recovery performance. In some works, the parallel between this problem and the coding/decoding paradigm is highlighted, and information theory tools are used to tackle the recovery, assuming a field structure for the alphabet \cite{das13, bio14sec}. In other works, CS techniques are rearranged for the finite-valued problem, for example, greedy pursuit algorithms \cite{bio14, spa15,fli18} or Bayesian methods \cite{tia09}. A particular focus is dedicated to binary signals, due also to their relevance in  digital communications problems \cite{tia09, sto10, spa15, lee16,fox18asi}.

{\textcolor{black}{ In the mentioned works, theoretical analyses are missing or limited to particular assumptions (\emph{e.g.}, the field structure).}}  This gap has been recently filled by \cite{kei17,fli18}. \cite{kei17} provides a CS theory for finite-valued signals: starting from the Basis Pursuit  formulation (BP, \cite{fou11}), the authors propose to impose the (convex) constraint $x\in \text{conv}(\alf^n)$ (where conv indicates the convex hull) and theoretically analyze the problem in terms of null space properties, robustness to noise, and phase transitions. {\textcolor{black}{ In \cite{fli18}, instead, a greedy pursuit algorithm called PROMP is developed and mathematically analyzed to recover sparse signals over lattices.}}
\subsection{Concave penalization}
At the same time, in (not finite-valued) CS and sparse signal processing, the use of \emph{concave} penalties has become increasingly popular in the last years, as it has been observed to be more efficient with respect to the convex formulation \cite{can08rew, gas09, woo16}. In many practical applications (such as X-ray CT  \cite{cha13} and MRI \cite{cha09}) concave penalties provide lower computation complexity and more accurate recovery, with less measurements with respect to the convex formulation. We notice that the enhancement thanks to concave penalties is not limited to linear acquisition problem, but has been observed in other machine learning frameworks: for example, in \cite{lap12}, a support vector machine with concave penalty is shown to produce a more parsimonious model without affecting the accuracy.
For what concerns CS, in \cite{woo16} some fundamental theoretical results are proven about BP and Basis Pursuit Denoising (BPDN) with concave penalty, in terms of properties of the minimum and robustness to noise.
\subsection{Our contribution}
The aim of this paper is to prove that the finite-valued sparse signal recovery problem (with special focus on CS) can be efficiently approached using a concave penalty. In particular, we obtain an enhancement with respect to the state-of-the-are \cite{kei17,fli18}. We provide both theoretical and experimental results. Our main contributions can be summarized in four points.

1)  \emph{Definition and analysis of a cost functional}: we define a non-convex Lasso-kind functional, formed by a least squares term plus a concave penalty, and we prove that its global minimum is the desired finite-valued sparse signal. As a difference from classical Lasso, our system is then unbiased. The analysis is based on  null space properties.

2) \emph{Robustness to noise}:  we analyze the proposed model in the presence of noise, both on the signal and on the measurements. We remark that we use the same  model for the noise-free and the noisy cases, which is suitable for systems where both exact and noisy measurements are expected.

3) \emph{Development of a recovery algorithm}: we derive an ADMM-based algorithm to search the desired minimum and we prove its convergence to the set of stationary points.

4) \emph{Validation via numerical experiments}: we conduct numerical experiments that highlight the efficiency of the proposed algorithm, in terms of estimation accuracy, {\textcolor{black}{required}} number of measurements, and convergence speed.

The paper is organized as follows. In Section \ref{sec:problem}, we define our model, and we put it into perspective with respect to prior literature. We start describing the bipolar ternary case $\alf=d\{0,\pm 1\}$, for some fixed $d>0$, and then extend to generic bipolar alphabets  $\alf=d\{0,\pm 1,\dots, \pm q\}$, $q\in\N$ (for simplicity, we assume that the symbols in $\alf$ are equidistant). In sections \ref{sec:ternary}-\ref{sec:generic} we prove theoretical results on the proposed model, for ternary and generic alphabets. In Section \ref{sec:algorithm}, we develop and analyze a recovery algorithm, which is further implemented and tested via numerical experiments in Section \ref{sec:numerical}. Finally, we draw some conclusions in the last section.
\subsection{Notation}
Throughout this paper, we use $I$ for the identity matrix (dimension is not specified when evident). $n$ and $S\subseteq \{1,\dots,n\}$ respectively are the length of the sparse signal and its support, while $k=|S|$ is the sparsity level (\emph{i.e.}, the number of non-zero components). $S^c$ is the complementary of $S$. For any $v\in\R^n$, $v_S$ is the restriction of $v$ on the components in $S$. $\R$ and $\R_+$ are the sets of real numbers and real non-negative numbers; $\Z$ and $\N$ are the sets of the integers and of the natural numbers. $B\succ 0$ means that $B$ is a positive definite matrix. Given a vector of weights $\zeta\in\R_+^n$, we define the weighted $\ell_1$-norm as $\|x\|_{1,\zeta}:=\sum_{i=1}^n \zeta_i |x_i|$.
\section{Problem statement}\label{sec:problem}
The problem of sparse signal recovery from linear measurements can be conceived as an $\ell_1$ convex minimization problem, known as Lasso \cite{tib96}: 
\begin{equation}\label{lasso}  
 \min_{x\in\R^n}\frac{1}{2}\left\|y-A x \right\|_2^2+\lambda \| x\|_1,~~~\lambda>0
\end{equation}
where $A\in\R^{m\times n}$ is the sensing matrix, and $\lambda>0$ is a parameter to set; in CS, $m<n$. The $\ell_1$-norm is known to well approximate the $\ell_0$-norm \cite{fou11} and has the important advantage of transforming the problem from combinatorial to convex, which makes it solvable in polynomial time. Iterative algorithms are often used to solve Lasso, \emph{e.g.}, Iterative Soft Thresholding (IST, \cite{dau04, for10}) and Alternating Direction Method of Multipliers (ADMM, \cite{boy10}). Both converge to the minimum of the Lasso functional. ADMM is known to require a significantly smaller number of iterations with respect to IST, keeping similar low complexity per iteration. This makes ADMM more attractive, along with its predisposition to distributed and parallel systems \cite{boy10, mata15, fia18}. 

We remark that the parameter $\lambda$ in Lasso has to be designed based on the noise: a larger $\lambda$ may tolerate a larger noise. Nevertheless, this has a drawback: Lasso has always a bias (proportional to $\lambda$) \cite{zha10MCP}, therefore the signal of interest is never exactly recovered in absence of noise. In classical sparse signal recovery and CS, the noise-free and noisy cases are then tackled using different models (in the convex setting, BP for the noise-free case, and BPDN or Lasso in the noisy case). This distinction is clearly not optimal for systems where acquisition is sometimes corrupted by noise and sometimes not. In this work, we instead show that a suitable concave penalization can remove this drawback for finite-valued signals.

\begin{figure}
\centering
 	\begin{tikzpicture}[scale=0.82]
		\begin{axis}[xmin=-1.6,xmax=2,ymin=-0.3, ymax=1.1, xtick={-1,0,1}, ytick={0,1}, ylabel={Penalty},axis lines=left, axis equal image, legend entries={\small{$\ell_0$}, ,\small{$\ell_1$}  ,\small{$\ell_p^p$}, \small{$\log$}, \small{MCP }},legend pos= south east]  
			\addplot[draw=red, mark=none, semithick, domain=-1.8:-0.05,-(,shorten >=-1pt] 
		{1};
			\addplot[ draw=red, mark=none, semithick, domain=1.8:0.05,-(,shorten >=-1pt] 
		{1};
		\addplot[name path = funone, draw=black, mark=none, semithick] {abs(x)};
		\addplot[draw=yellow, mark=none, thick,samples=5000] {sqrt(abs(x))};
 		\addplot[draw=green, mark=none, semithick,samples=1500] {0.2*ln( 100*abs(x)+1)}; 
 		\addplot[draw=blue, mark=none,thick,domain=-sqrt(2):sqrt(2)]{abs(x)*sqrt(2)-0.5*x^2} ; 
 		\addplot[ draw=blue, mark=none, thick, domain=-1.8:-sqrt(2)] {1};
 		\addplot[ draw=blue, mark=none, thick, domain=1.8:sqrt(2)] {1}	;
 		\draw[color=red, fill] (axis cs:0,0) circle (2.3pt);
\end{axis}
	\end{tikzpicture}
	\caption{Some popular concave penalties. Concave penalties are closer to $\ell_0$ than $\ell_1$.} 
	\label{fig:penalty}
\end{figure}
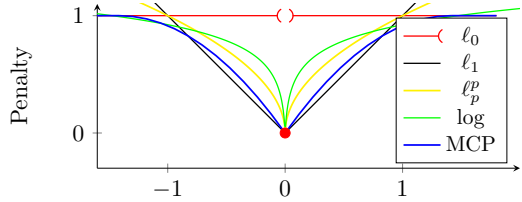

We now reformulate Lasso using a concave penalty \cite{woo16}:
\begin{equation}\label{concave_penalization}  
\begin{split}
 &\min_{x\in\R^n}\frac{1}{2}\left\|y-A x \right\|_2^2+\lambda \sum_{i=1}^n g(|x_i|)\\&g:\R_+\to \R_+ \text{ concave, nondecreasing in } |x_i|.
\end{split}
\end{equation}
The intuition behind the success of concave penalization is that concave functions approximate the $\ell_0$-norm  better than $\ell_1$, as one can see in Figure \ref{fig:penalty}. In the overdetermined setting, concave penalization has been largely studied within the statistical community \cite{fan01_pioneer, fan11, fan14, zou08_LLA, zha10MCP, zha12, gas09}. In these papers, conditions to have the oracle property and to reduce the Lasso bias are studied, mainly in the asymptotic  case $n\to\infty$ \cite{fan01_pioneer, fan11}. In the concave CS literature, instead, the following BP/BPDN formulations are commonly used:
\begin{equation}\label{eq:concave_bp}
\begin{split}
&\text{Concave BP: }\min \sum_{i=1}^n g(|x_i|)\text{ s.t. }Ax=y\\
&\text{Concave BPDN: }\min \sum_{i=1}^n g(|x_i|)\text{ s.t. } \left\|Ax-y\right\|<\eta,~\eta>0\\
&g:\R_+\to \R_+ \text{ concave, nondecreasing in } |x_i|.
\end{split}
\end{equation}
The most popular $g$'s belong to these families:
\begin{itemize}
	\item $\ell_p$, with $p\in (0,1)$ \cite{irls, rav15irls, woo16};
	\item log-sum: $\log |x_i|+\varepsilon$, $\varepsilon>0$ \cite{can08rew, faz03, cal16};
	\item smoothly clipped absolute deviations (SCAD) \cite{fan01_pioneer}; 
	\item minimax concave penalty (MCP) \cite{zha10MCP, woo16, fox16, hua18}.
\end{itemize}
We remark that in some works the concave penalization arises from iterative re-weighting strategies \cite{can08, irls, rav15irls}, which have been largely studied in the last years. Others works, instead, consider algorithms based on difference of convex (DC) functions programming \cite{gas09}. Both re-weighting and DC approaches, however, imply the solution of a convex problem (typically a Lasso) during each iteration, which makes them more complex than the proposed approach (see Section \ref{sec:algorithm}).

In this paper, {\textcolor{black}{we focus on MCP and we show that it is an optimal choice to deal with finite-valued signals}}. The literature on MCP in sparse signal recovery is widespread. In \cite{woo16}, concave BP/BDN \eqref{eq:concave_bp} with MCP is studied and shown to perform better than  $\ell_p$, $p\in (0,1)$. In  \cite{cha09, cha13}, numerical experiments support its good performance. In \cite{fox16}, MCP is shown to well adapt to a distributed recovery setting. In \cite{zha10MCP}, MCP is exploited to build an Lasso-kind estimator with reduced bias. In  machine learning, \cite{lap12} uses MCP  to build a concave support vector machine for parsimonious feature selection.

Our use of MCP is now described, starting from the ternary alphabet $\{0,\pm d\}$, $d>0$, and then extending to generic alphabets.
\subsection{Ternary finite-valued signals}
Let $\xtrue\in\{0,\pm d\}^n$, where $d>0$ is known. MCP is defined as follows: for any $z\in\R$,
\begin{equation}\label{my_g}
g(|z|):=\left\{\begin{array}{lr}
d |z|-\frac{1}{2}z^2&\text{ if } z\in[-d,d]\\
\frac{1}{2}d^2 &\text{ otherwise.}\\
\end{array}\right.
\end{equation}
Let us consider the following cost functional in the convex hull $[-d,d]^n$:
\begin{equation}\label{mcplasso}  
\begin{split}
 \fun(x)&:=\frac{1}{2}\left\|y-A x \right\|_2^2+\lambda d \left\|x\right\|_1-\frac{\lambda}{2}\left\|x \right\|_2^2\\
 &\lambda>0,~x\in[-d,d]^n.
\end{split}
\end{equation}
In the following, we sometimes use the notation $G(x):=d \left\|x\right\|_1-\frac{1}{2}\left\|x \right\|_2^2$, and we  refer to $\fun(x)$ as  MCP-Lasso.

{\textcolor{black}{An intuitive motivation to expect that model \eqref{mcplasso} is suitable for signals in $\{0,\pm d\}^m$ is the presence of the term $-\left\|x \right\|_2^2$: given a sparse signal $v\in[-d,d]^n$ with fixed $\|v\|_1$, its energy is maximized (that is,  $-\left\|v \right\|_2^2$ is minimized) by pushing the non-zero entries to the boundaries $-d$ and $d$. This intuition is rigorously proven by our theoretical analysis in next sections.}}
\subsection{Generic finite-valued signals}
Let $\xtrue \in \mathcal{A}^m$ where $\mathcal{A}=d\{0,\pm 1, \pm 2,\dots, \pm q\}$. We propose to reformulate MCP-Lasso as follows:
\begin{equation}\label{mcplasso_gen}  
\begin{split}
 &\gun(x):=\frac{1}{2}\left\|y-A x \right\|_2^2+\lambda \sum_{i=1}^n \beta_i(x_i)~|x_i|-\lambda\frac{1}{2}\left\|x \right\|_2^2\\
&~\lambda>0,~x\in [-dq,dq]^n,~\beta_i(x_i):=\min\{\alpha\in\mathcal{A}~ \text{s.t. } |x_i|\leq \alpha \}.
\end{split}
\end{equation}
{\textcolor{black}{In this formulation, each weight $\beta_i(x_i)$ increases when  the magnitude of $x_i$ increases, according to the generic principle of reweighting methods. Moreover, we remark that the given definition of the $\beta_i(x_i)$'s preserves the non-negativity of the penalization. This might preserved even choosing, for example, all the weights equal to $qd$; however, even though not intuitive, the proposed quantization is the key to characterize the minimum of the functional and its suitability for the signal recovery, as will be clear in Section \ref{sec:generic}.}}

%
{\textcolor{black}{\begin{remark}
In this paper, only bipolar alphabets are considered, which have both positive and negative symbols \cite{kei17}. However, as observed in \cite{kei17,sto10}, the unipolar case, where symbols all have the same sign, presents specific features that improve the performance in the BP approach. The MCP approach restricted to the alphabet $\{0,1\}$ is partially investigated in \cite{fox18asi}, where a characterization of local minima is provided. No particular improvement has been yet highlighted for unipolar alphabets using MCP. However, this topic is currently under analysis.   
\end{remark}}}
{\textcolor{black}{We remark that $\fun$ and $\gun$ are in general non-convex, hence minimization might be complicated. However, they have a semi-algebraic expression, therefore global minimization can be performed via semi-algebraic optimization methods \cite{las15}, which marks an advantage with respect to other concave penalties different from MCP. Since these methods are numerically complex in the large scale, in this paper we develop iterative strategies for minimization (see Section \ref{sec:algorithm}).}}
\begin{remark}\label{marginally}
If $A^T A -\lambda I \succ 0$, $\fun$ and $\gun$ are convex. However, this condition does not match with CS, which is our main focus. We then discuss it only marginally.
\end{remark}
\subsection{Relation to prior literature}
In this paper, the main works that we refer to are \cite{woo16} for what concerns concave penalization, and \cite{kei17} for what concerns finite-valued sparse signal recovery.

As already mentioned, in \cite{woo16}, BP and BPDN  are recast to \eqref{eq:concave_bp} using $\ell_p$, $p\in(0,1)$, and MCP. The authors prove that (a) under technical conditions, the true signal is the minimum of concave BP in the noise-free case; (b) concave BPDN is robust to noise. Their analysis is mainly based on null space properties  \cite[Chapter 4]{fou11}. They also derive an IST procedure, which is convergent for  $\ell_p$ penalty, under the condition $A^T A -\lambda I \succ 0$. The problems of the convergence in the MCP case and in the CS setting are then open. No numerical results are shown in \cite{woo16}.

In \cite{kei17}, instead, the problem of CS with  finite-valued sparse signals is theoretically analyzed, and recast to a convex formulation by considering the convex hull of the alphabet. Null space properties, robustness, and phase transitions are discussed. Some numerical results are shown, where recovery is performed via convex programming.

As in \cite{woo16, kei17}, our analysis is based on null space conditions and robustness bound are provided. Similarly to \cite{woo16}, we prove that in the noise-free setting the global minimum is the true signal. As a difference from \cite{woo16}, we use a unique model/functional for the noise-free and noisy settings, and we provide an efficient recovery algorithm. On the other hand, we improve the performance with respect to \cite{kei17} by exploiting the concavity.
\section{Theoretical analysis: ternary alphabet}\label{sec:ternary}
In this section, we characterize the minimum of $\fun(x)$ defined in \eqref{mcplasso} and we prove its robustness to noise. Let us state a preliminary lemma.
\begin{lemma}\label{doublereg}
Let $\xtrue\in\{0,\pm d\}^n$ be $k$-sparse with support $S$, and let $y=A\xtrue$. Let $\xmin$ be the global minimum of $\fun(x)$ defined in \eqref{mcplasso} over $[-d,d]^n$. Then,
\begin{equation}
\begin{split}
&\text{(a)}~~\left\|\xmin\right\|_2 \leq \left\|\xtrue\right\|_2;~~~~\text{(b)}~~\left\|\xmin\right\|_1 \leq \left\|\xtrue\right\|_1.\\
\end{split}
\end{equation}
\end{lemma}
\begin{proof}
Let us notice that 
\begin{equation}\label{ternary_property}
G(\xtrue)=\frac{1}{2}kd^2=\frac{1}{2}\left\|\xtrue\right\|_2^2=\frac{1}{2}d\left\|\xtrue\right\|_1=\frac{1}{\lambda} \fun(\xtrue).
\end{equation}
Since $\fun(\xmin)\leq \fun(\xtrue)$ and $\lambda G(z)\leq \fun(z)$ for any $z\in\R^n$, then $\lambda G(\xmin)\leq \fun(\xmin) \leq \fun(\xtrue) = \lambda G(\xtrue)$. {\textcolor{black}{Now, let us notice that for any $x\in[-d,d]^n$, $G(x)\geq \frac{1}{2}\left\|x\right\|_2^2$: since $G(x)=d\|x\|_1-\frac{1}{2}\|x\|_2^2$, the inequality is true if $d\|x\|_1\leq \|x\|_2^2$, which in turn is true because for each $x_i\in[-d,d]$, $d|x_i|\geq x_i^2$.}}

We apply this to state that $G(\xmin)\geq \frac{1}{2}\left\|\xmin\right\|_2^2$, which in turn can be substituted into the previous inequality $G(\xmin)\leq G(\xtrue)$ to finally obtain (a). To obtain (b), it suffices to write explicitly $G(\xmin)\leq G(\xtrue)$ and apply (a).
\end{proof}
\subsection{Global minimum of $\fun$}
We now prove that the minimum of $\fun(x)$ corresponds to the desired signal in the ternary case under null space conditions, commonly used in CS \cite[Chapter 4]{fou13}. Let us recall the following definition.
 \begin{definition}\label{RNSP}\cite[Definition 4.21]{fou13}
A matrix $A\in\R^{m,n}$ is said to satisfy the $\ell_2$-robust null space property (\emph{RNSP}) of order $k$ with parameters $\rho\in(0,1)$, $\tau>0$ if
$$\left\|v_S\right\|_2 \leq \frac{\rho}{\sqrt{k}} \left\|v_{S^c}\right\|_1+\tau \left\|Av\right\|_2$$ for any $v\in\R^n$ and for any $S\subset \{1,\dots,n\}$ with $|S|\leq k$.
\end{definition}
\begin{theorem}\label{theo:globalminimum}
Let us consider $y=A\xtrue$, where $\xtrue$ is a $k$-sparse signal in $\{0,\pm d\}^n$, and $A\in\R^{m,n}$ with $m<n$.  If $A$  satisfies the \emph{RNSP} of order $k$ with parameters  $\tau>0$ and $\rho\in\left(0,\sqrt{1-\lambda\tau^2}\right)$ (provided that $\lambda<\tau^{-2}$), then 
$\xtrue$ is the {\textcolor{black}{unique}} global minimum of $\fun$ defined in \eqref{mcplasso} over $[-d,d]^n$.
\end{theorem}
{\textcolor{black}{We remark that the requirements $\rho\in\left(0,\sqrt{1-\lambda\tau^2}\right)$ and $\lambda<\tau^{-2}$ are not restrictive, as $\lambda$ is a design parameter that can be set as small as necessary in the noise-free case.}}
\begin{proof}
Let $\xtrue_i\in\{0,\pm d\}$ with support $S$. Consider any $h_i\in\R$ such that  $\xtrue_i+h_i\in [-d,d]$.
 \begin{equation*}
  \begin{split}
   &\fun(\xtrue+h)  = \frac{1}{2}\left\|A h \right\|_2^2+\lambda d \left\|\xtrue +h\right\|_1-\lambda\frac{1}{2}\left\|\xtrue +h\right\|_2^2\\
    & = \frac{1}{2}\left\|A h \right\|_2^2+\lambda d \left\|\xtrue +h\right\|_1 -\frac{\lambda}{2}\left\|\xtrue\right\|_2^2-\frac{\lambda}{2}\left\|h\right\|_2^2-\lambda\langle \xtrue,h \rangle.
  \end{split}
 \end{equation*}
 Let us focus on the term $ d \left\|\xtrue +h\right\|_1-\langle \xtrue,h \rangle=\sum_{i=1}^n  d |\xtrue_i+h_i|-\xtrue_i h_i$. We distinguish the following cases:
 \begin{itemize}
  \item if $\xtrue_i= d$ and $h_i\in[- d,0]$, then $ d |\xtrue_i+h_i|-\xtrue_i h_i = d^2$;
  \item if $\xtrue_i= d$ and $h_i\in[-2 d,- d]$, then $ d |\xtrue_i+h_i|-\xtrue_i h_i =- d^2-2 d h_i\geq  d^2$;
  \item if $\xtrue_i=- d$ and $h_i\in[0, d]$, then $ d |\xtrue_i+h_i|-\xtrue_i h_i = d^2$;
  \item if $\xtrue_i=- d$ and $h_i\in[ d,2 d]$, then $ d |\xtrue_i+h_i|-\xtrue_i h_i =- d^2+2 d h_i\geq  d^2$;
  \item if $\xtrue_i=0$ and $h_i\in[- d, d]$, then $ d |\xtrue_i+h_i|-\xtrue_i h_i = d|h_i|$.
 \end{itemize}
Therefore, 
\begin{equation*}
  \begin{split}
   \fun(\xtrue+h) & \geq \frac{1}{2}\left\|A h \right\|_2^2-\frac{\lambda}{2}\left\|\xtrue\right\|_2^2-\frac{\lambda}{2}\left\|h\right\|_2^2\hspace{-0.1cm}+\hspace{-0.1cm}\lambda k  d^2\hspace{-0.1cm} +\hspace{-0.1cm} \lambda d\left\|h_{S^c}\right\|_1 \\
   & =\fun(\xtrue) +\frac{1}{2}\left\|A h \right\|_2^2-\frac{\lambda}{2}\left\|h\right\|_2^2+\lambda  d \left\|h_{S^c}\right\|_1 \\
  \end{split}
 \end{equation*}
 where we use the fact that $\fun(\xtrue)=\frac{\lambda}{2}\left\|\xtrue\right\|_2^2=\frac{\lambda}{2}k d^2$. In order to prove the thesis, it is then sufficient to prove that, for any $h\in[-d,d]^n\setminus\{0\}$,
\begin{equation}\label{ineedpos}
\frac{1}{2}\left\|A h \right\|_2^2-\frac{\lambda}{2}\left\|h\right\|_2^2+\lambda  d \left\|h_{S^c}\right\|_1>0. 
\end{equation}
For any $h\in[-d,d]^n$, it is straightforward to prove the following inequality:
\begin{equation}\label{usefulobs}
d \left\|h_{S^c}\right\|_1-\frac{1}{2}\left\|h_{S^c}\right\|_2^2\geq \frac{d}{2} \left\|h_{S^c}\right\|_1.
\end{equation}
{\textcolor{black}{We now distinguish two cases.}}

{\textcolor{black}{\textbf{Case 1:} $\|h_{S^c}\|_1> kd$. Plugging \eqref{usefulobs} into the left-hand expression in \eqref{ineedpos}, we obtain $ \frac{1}{2}\left\|A h \right\|_2^2+ \frac{\lambda d}{2} \left\|h_{S^c}\right\|_1-\frac{\lambda}{2}\left\|h_S\right\|_2^2$. Since $\left\|h_S\right\|_2^2\leq kd^2$, we have $d \left\|h_{S^c}\right\|_1-\left\|h_S\right\|_2^2>kd^2-kd^2=0$, which implies \eqref{ineedpos}.}}

{\textcolor{black}{\textbf{Case 2:} $\|h_{S^c}\|_1\leq kd$. We exploit \eqref{usefulobs} and the RNSP.
\begin{equation*}
{\small{\begin{split}
&\frac{1}{\lambda}\left\|A h \right\|_2^2-\left\|h\right\|_2^2+2 d \left\|h_{S^c}\right\|_1\geq \frac{1}{\lambda}\left\|A h \right\|_2^2+d\left\|h_{S^c}\right\|_1-\left\|h_S\right\|_2^2\\
&\geq\frac{1}{\lambda}\left\|A h \right\|_2^2+d \left\|h_{S^c}\right\|_1-\left[ \frac{\rho}{\sqrt{k}} \left\|h_{S^c}\right\|_1+\tau \left\|Ah\right\|_2 \right]^2\\
\end{split}}}
\end{equation*}
Since $\|h_{S^c}\|_1^2\leq kd\|h_{S^c}\|_1$, the last expression is not smaller than
$$\left(\frac{1}{\lambda}-\tau^2\right)\left\|A h \right\|_2^2+d\left(1-\rho^2  \right) \left\|h_{S^c}\right\|_1-\frac{2\rho\tau}{\sqrt{k}} \left\|h_{S^c}\right\|_1 \left\|Ah\right\|_2.$$
For simplicity, let us name $c_1=\left(\lambda^{-1}-\tau^2\right)$, $c_2=d\left(1-\rho^2  \right)$, and $c_3=\frac{2\rho\tau}{\sqrt{k}}$. The last expression is then equal to
\begin{equation*}
\begin{split}
&c_1\left\|A h \right\|_2^2+c_2 \left\|h_{S^c}\right\|_1-c_3 \left\|h_{S^c}\right\|_1 \left\|Ah\right\|_2\pm \frac{c_3^2}{4c_1}\left\|h_{S^c}\right\|_1^2\\
&=\left(\sqrt{c_1}\left\|A h \right\|_2-\frac{c_3}{2\sqrt{c_1}}\left\|h_{S^c}\right\|_1\right)^2+\left(c_2-\frac{c_3^2}{4c_1}kd\right)\left\|h_{S^c}\right\|_1\\
&\geq\left(c_2-\frac{c_3^2}{4c_1}kd\right)\left\|h_{S^c}\right\|_1^2= d\left(1-\frac{\rho^2}{1-\lambda\tau^2}\right)\left\|h_{S^c}\right\|_1\\
\end{split}
\end{equation*}
which is positive for $\rho\in(0,\sqrt{1-\lambda\tau^2})$ (where $\lambda<\tau^{-2}$). This concludes the proof.}}
\end{proof}
{\textcolor{black}{\begin{remark}\label{qf}
Even though very popular in CS, the RNSP is not easy to check for a given $A$. In the literature, its validity is proven for families of matrices, like Gaussian matrices \cite[Chapter 9]{fou13}. This limits the practical use of results based on RNSP. For this motivation, we notice that $d \left\|h_{S^c}\right\|_1-\frac{1}{2}\left\|h_{S^c}\right\|_2^2\geq \frac{d}{2} \left\|h_{S^c}\right\|_2^2$ for any $h\in[-d,d]^n$, thus
$$\frac{1}{\lambda}\left\|A h \right\|_2^2-\left\|h\right\|_2^2+2 d \left\|h_{S^c}\right\|_1\geq \frac{1}{\lambda}\left\|A h \right\|_2^2+d \left\|h_{S^c}\right\|_2^2-\left\|h_S\right\|_2^2.$$
To prove Theorem \ref{theo:globalminimum}, it is then sufficient to check whether the quadratic form $\lambda^{-1}\left\|A h \right\|_2^2+d \left\|h_{S^c}\right\|_2^2-\left\|h_S\right\|_2^2$ is positive, that is, whether the eigenvalues of $\lambda^{-1}A^T A+d I_{S^c}-I_S$ are positive, where $I_S, I_{S^c}\in\{0,1\}^{n,n}$  have entries equal to one on the diagonal on the positions corresponding to $S$ and $S^c$, respectively (and zero otherwise). As $S$ is not priorly known, one has to check all the possible supports, which of course might take long time. However, this is definitely more feasible than checking the RNSP, which requires to solve a non-convex semi-algebraic optimization problem.
\end{remark}
\begin{remark}
In CS, it is known that the RNSP is related to the restricted isometry property \cite[Theorem 6.13]{fou13}, on the basis of which one can estimate that the number of necessary measurements  $m$ to have the RNSP is (asymptotically) of order $k\log(en/k)$ for random Gaussian and Bernoulli matrices \cite[Theorem 9.27]{fou13}.
\end{remark}}}
{\textcolor{black}{Even though this work is focused on the compressed case, it is worth to mention the following result for the non-compressed case, which can easily obtained following the proof of Theorem \ref{theo:globalminimum} (see also Remark \ref{marginally}).}}
\begin{corollary}\label{cor:noncs}
If $A^T A -\lambda I \succ 0$, $\fun$ is convex and $\xtrue$ is its global minimum.
\end{corollary}
\begin{remark}
Theorem \ref{theo:globalminimum} and Corollary \ref{cor:noncs} highlight a substantial difference between Lasso and MCP-Lasso. In the noise-free case, Lasso is biased: the true signal never corresponds to the minimum, with a bias proportional to $\lambda$. Indeed, Lasso is conceived for the noisy case: a larger $\lambda$ may tolerate a larger noise. In classical CS, noise-free and noisy cases are approached separately, with different models (\emph{e.g.}, BP vs BPDN/Lasso). This is not optimal  for systems where both noise-free and noisy measurements are acquired. This problem is overcome by MCP-Lasso, where $\lambda$ could be designed on the maximum noise, without bias when noise-free measurements are acquired.
\end{remark}
\subsection{Robustness to noise of $\fun$}
We know analyze the robustness to noise of $\fun(x)$ under the  RNSP. We consider 
$$y=A(\xtrue+\delta)+\epsilon$$ where $\delta\in\R^n$ and $\epsilon\in\R^m$ respectively represent the signal noise (\emph{i.e.}, the signal is not exactly sparse) and the measurement noise.  We now prove a robust bound of kind \cite[Theorems 4.19, 4.21]{fou13} for the distance between the desired $\xtrue$ and the global minimum of $\fun$.

Rearranging Lemma \ref{doublereg} in the noisy setting, we can prove the following inequalities (we omit the proof, which can be simply derived from the proof of Lemma \ref{doublereg}).
\begin{lemma}\label{doublereg_noise}
\begin{equation*}
\begin{split}
&(a)~~~\left\| \xmin \right\|_1 \leq \left\| \xtrue\right\|_1+ \frac{1}{2\lambda}\left\| A\xtrue-y \right\|_2^2-\frac{1}{2\lambda}\left\| A\xmin-y \right\|_2^2;\\
&(b)~~~\left\| \xmin \right\|_2^2 \leq \left\| \xtrue\right\|_2^2+ \frac{1}{\lambda}\left\| A\xtrue-y \right\|_2^2-\frac{1}{\lambda}\left\| A\xmin-y \right\|_2^2.
\end{split}
\end{equation*}
\end{lemma}

\begin{theorem}\label{robustheorem}
Let $\xtrue\in\{0\pm d\}^n$ be a $k$-sparse signal with support $S$, and let $\xtrue+\delta$ be its noisy version (i.e., the signal is not exactly sparse). Let $y=A(\xtrue+\delta)+\epsilon$, where $\epsilon$ is the measurement noise.  If $A\in\R^{m,n}$ satisfies the \emph{RNSP} of order $k$ with parameters $\rho\in(0,1), \tau>0$ , then,
{\small{$$\left\| \xmin - \xtrue\right\|_1 \leq \frac{1+\rho}{(1-\rho)2\lambda}\left\| A\delta+\epsilon\right\|_2^2+\frac{4\tau}{1-\rho}\left[\left\| A\delta+\epsilon\right\|_2+ \frac{d}{2}\sqrt{k\lambda}\right]\hspace{-0.1cm}.$$}}
\end{theorem}
\begin{proof}
Let $h=\xmin - \xtrue.$ By triangle inequality, we have:
\begin{equation}\label{wtri}
\left\| \xtrue_S\right\|_1\leq \left\| \xtrue_S +h_S\right\|_1 + \left\| h_S\right\|_1.
\end{equation}
Moreover, it is easy to verify that we can decouple as follows:
\begin{equation}\label{wdecouple}
\left\| \xmin\right\|_1 = \left\| \xmin_S \right\|_1+\left\| \xmin_{S^c} \right\|_1= \left\| \xmin_S \right\|_1+\left\| h_{S^c} \right\|_1
\end{equation}
as $\xtrue_{S^c}=0$. Keeping in mind \eqref{wtri} and \eqref{wdecouple}, we can write:
\begin{equation*}
\begin{split}
\left\| h_{S^c} \right\|_1 & = \left\| \xmin \right\|_1 -\left\| \xmin_S \right\|_1 = \left\| \xmin \right\|_1 -\left\| \xtrue_S+h_S \right\|_1\\
& \leq \left\| \xmin  \right\|_1 - \left\| \xtrue\right\|_1+\left\|h_S \right\|_1.\\
\end{split}
\end{equation*}
Exploiting the RNSP,
\begin{equation*}
\begin{split}
&\left\| h_{S^c} \right\|_1 \leq  \left\| \xmin  \right\|_1 - \left\| \xtrue\right\|_1+\rho\left\|h_{S^c} \right\|_1+ \tau\left\|Ah \right\|_2\\
\Rightarrow & \left\| h_{S^c} \right\|_1 \leq \frac{1}{1-\rho}\left[ \left\| \xmin  \right\|_1 - \left\| \xtrue\right\|_1+ \tau\left\|Ah \right\|_2\right].
\end{split}
\end{equation*}
Therefore,
\begin{equation}\label{qui}
\begin{split}
\left\| h \right\|_1 &= \left\| h_S \right\|_1 +\left\| h_{S^c} \right\|_1 \leq (1+\rho)\left\| h_{S^c} \right\|_1 + \tau \left\|A h\right\|_2\\
&\leq \frac{1+\rho}{1-\rho}\left[ \left\| \xmin  \right\|_1 - \left\| \xtrue\right\|_1\right]+ \frac{2}{1-\rho}\tau\left\|Ah \right\|_2.
\end{split}
\end{equation}
According to Lemma \ref{doublereg_noise}.(a), $\left\| \xmin  \right\|_1 - \left\| \xtrue\right\|_1\leq \frac{1}{2\lambda}\left\| A\xtrue-y \right\|_2^2$. Moreover, by Lemma \ref{doublereg_noise}.(b):
{\small{\begin{equation}\label{eq:Ah}
\left\|Ah \right\|_2 \leq \left\|A\xtrue-y \right\|_2+ \left\|A\xmin-y \right\|_2 \leq 2\left\|A\xtrue-y \right\|_2+ \sqrt{\lambda}\left\| \xtrue\right\|_2.
\end{equation}}}
Plugging the last two inequalities  into \eqref{qui} we conclude our proof.
\end{proof}
\begin{remark}\label{etalambda}
In line with classical robustness results \cite{fou13, woo16}, Theorem \ref{robustheorem} states that the recovery error is driven by the noises $\delta$ and $\epsilon$, and by the parameter $\lambda$, which is designed based on the noise magnitude.
For example, if $\left\| A\delta+\epsilon\right\|_2\leq \eta $ for some known $\eta>0$, one can set $\lambda=\eta$ and the bound becomes $\left\| h \right\|_1\leq \eta \left[\frac{1+\rho}{2(1-\rho)}+\frac{2\tau}{1-\rho}\right]+ d\sqrt{k\eta} d.$ In this way, when the noise tends to zero, then also the recovery error tends to zero.
 \end{remark}
\begin{remark}
In \cite{woo16} two robustness results are proven for the (non-finite valued) concave problem. First, \cite[Proposition 4.4]{woo16} requires the so-called G-NNSP \cite[Definition 4.3]{woo16}, which might be more stringent. Second, \cite[Theorem 4.5]{woo16} provides a noise-driven bound when $AA^T=I$ assuming that $k g(2 \beta')< (n-k)g(\alpha')$ where $\alpha'$ and $\beta'$ respectively are the minimum and the maximum magnitudes of the projection of $h$ onto ker$(A)$ ($\alpha'>0$ is guaranteed only if $\delta>0$ \cite[Lemma 4.2]{woo16}, that is, the signal must be not exactly sparse). Since $\alpha'$ and $\beta'$ depend on $h$, which is the quantity to be estimated, this result is somehow circular.
\end{remark}
\section{Theoretical analysis: generic alphabet}\label{sec:generic}
In this section, we extend the theoretical results to a generic bipolar alphabet $\alf=d\{0,\pm 1,\dots, \pm q\}$, $d>0$, $q\in\N$.

\subsection{Global minimum of $\gun$}
\begin{theorem}\label{theo:globalminimum_gen}
Let us consider $y=A\xtrue$, where $\xtrue$ is a $k$-sparse signal in $\alf^n$,  $\alf=d\{0,\pm 1,\dots, \pm q\}$, and $A\in\R^{m,n}$ with $m<n$. 
 If $A\in\R^{m,n}$ satisfies the RNSP of order $k$ with parameters $\tau>0$ and $\rho\in\left(0,\frac{\sqrt{1-2q\lambda\tau^2}}{2q^2}\right)$, and given $\lambda<\frac{1}{2q \tau^2}$, then $\xtrue$ is the {\textcolor{black}{unique}} global minimum of $\gun(x)$ defined in \eqref{mcplasso_gen} over $conv(\mathcal{A}^n)$.
\end{theorem}
\begin{proof}
Let $\xtrue_i\in\mathcal{A}$. Assume that  $\xtrue_i+h_i\in conv(\mathcal{A})$. Let $S$ be the support of $\xtrue$. 
Following the procedure of the proof of Theorem \ref{theo:globalminimum}, we prove that $\gun(\xtrue+h) -\gun(\xtrue)>0$ for any $h\in[-qd,qd]^n\setminus\{0\}$. 
 \begin{equation*}
 \begin{split}
 \gun(\xtrue+h)  &=  \gun(\xtrue)+\frac{1}{2}\left\|A h \right\|_2^2-\frac{\lambda}{2}\left\|h \right\|_2^2+\lambda \sum_{j\in S^c}\beta_j(h_j)|h_j|\\&+\lambda \sum_{i\in S}\left[\beta_i(\xtrue_i +h_i)~ |\xtrue_i +h_i| -\xtrue_i^2-\xtrue_i h_i \right].
   \end{split}
 \end{equation*}
 Let us focus on the terms $r_i=\beta_i(\xtrue_i +h_i)~ |\xtrue_i +h_i| -\xtrue_i^2-\xtrue_i h_i$, $i\in S$. We distinguish the following cases:
 \begin{itemize}
  \item if $\xtrue_i>0$, $h_i>-d$: $\beta_i(\xtrue_i +h_i)\geq\xtrue_i$, then $r_i> 0$;
  \item if $\xtrue_i>0$, $h_i\in[-\xtrue,-d)$, $r_i= (\beta_i(\xtrue_i +h_i)-\xtrue_i)(\xtrue_i+h_i) \geq h_i(\xtrue+h_i)=h_i^2+\xtrue_i h_i\geq h_i^2+q d h_i\geq(1-q)h_i^2$;
  \item if $\xtrue_i>0$, $h_i\in [-2\xtrue,-\xtrue)$: $\xtrue_i\geq\beta_i(\xtrue_i +h_i)$, then $r_i= -(\xtrue_i-\beta_i(\xtrue_i +h_i))(\xtrue_i+h_i) \geq 0$;
\item if $\xtrue_i>0$, $h_i<-2\xtrue$: $\beta_i(\xtrue_i +h_i)>\xtrue_i$, then $r_i\geq 0$.
 \end{itemize}
We omit the description of the case $\xtrue_i<0$ which is symmetric. In conclusion, $\gun(\xtrue+h)-\gun(\xtrue)$ is not smaller than
{\begin{equation*}
\begin{split}    
   &\frac{1}{2}\left\|A h \right\|_2^2-\frac{\lambda}{2}\left\|h \right\|_2^2+\lambda \sum_{j\in S^c}\beta_j(h_j)|h_j|-\lambda (q-1)\left\|h_S \right\|_2^2\\
   &\geq \frac{1}{2}\left\|A h \right\|_2^2+\frac{\lambda}{2} \sum_{j\in S^c}\beta_j(h_j)|h_j|-\lambda q\left\|h_S \right\|_2^2.
   \end{split}
 \end{equation*}
We now distinguish two cases.

\textbf{Case 1:} $\left\|h_{S^c}\right\|_1>2 q^3kd$. Since $\beta_j(h_j)|h_j|\geq d|h_j|$, we simply obtain: $\sum_{j\in S^c}\beta_j(h_j)|h_j|-2 q\left\|h_S \right\|_2^2>2q^3kd^2-2 q\left\|h_S \right\|_2^2\geq 2q^3kd^2-2qk(qd)^2=0.$

\textbf{Case 2:}  $\left\|h_{S^c}\right\|_1\leq2 q^3kd$}. The procedure is analogous to the Case 2 in the proof of Theorem \ref{theo:globalminimum}.
\begin{equation*}
\begin{split}
&\frac{1}{\lambda}\left\|A h \right\|_2^2+\sum_{j\in S^c}\beta_j(h_j)|h_j|-2 q\left\|h_S \right\|_2^2\\
&\geq\frac{\left\|A h \right\|_2^2}{\lambda}+ \sum_{j\in S^c}\beta_j(h_j)|h_j|-2 q\left[ \frac{\rho}{\sqrt{k}} \left\|h_{S^c}\right\|_1+\tau \left\|Ah\right\|_2 \right]^2
\end{split}
\end{equation*}
\begin{equation}\label{qui}
\begin{split}
&\geq c_1 \left\|A h \right\|_2^2+c_2 \left\|h_{S^c}\right\|_1 - c_3 \left\|h_{S^c}\right\|_1 \left\|Ah\right\|_2\\
&\geq\left(c_2-\frac{c_3^2}{4c_1}2q^3kd\right)\left\|h_{S^c}\right\|_1
\end{split}
\end{equation}
where $c_1=\lambda^{-1} -2q\tau^2$, $c_2=d \left(1-4q^4\rho^2  \right)$, and $c_3=2q\lambda\rho\tau/\sqrt{k}$.
As $c_2-\frac{c_3^2}{4c_1}2q^3kd=  d\left(1-\frac{4q^4\rho^2}{1-2q\lambda\tau^2}\right)$, the last expression in \eqref{qui} is positive when $\rho\in\left(0,\frac{\sqrt{1-2q\lambda\tau^2}}{2q^2}\right)$, provided that $\lambda<\frac{1}{2q \tau^2}$.
\end{proof}
\begin{remark}\label{qf2}
Similarly to Remark \ref{qf}, we notice that  $\lambda^{-1}\left\|A h \right\|_2^2\left\|h \right\|_2^2+2\lambda \sum_{j\in S^c}\beta_j(h_j)|h_j|-2 (q-1)\left\|h_S \right\|_2^2\geq \lambda^{-1}\left\|A h \right\|_2^2+\left\|h_{S^c} \right\|_2^2-2\left(q+\frac{1}{2}\right) \left\|h_S \right\|_2^2$. Therefore, the RNSP can be substituted by the study of the sign of a quadratic form, and, specifically, with the computation of the eigenvalues of $\lambda^{-1}A^T A+I_{S^c}-q I_S$. This can be useful in the practice, since the RNSP is difficult to prove for a matrix.
\end{remark}

\subsection{Robustness to noise of $\gun$}
Let us extend Lemma \ref{doublereg_noise} to the case of generic alphabet. The proof follows the schem of the proof of Lemma  \ref{doublereg_noise}, then omitted for brevity.
\begin{lemma}\label{doublereg_gen_noise}
Let $\xmin$ be the global minimum of $\gun$ over $[-qd,qd]^n$ with $y=A(\xtrue+\delta)+\epsilon$, $\xtrue\in\alf^n$. Then
\begin{equation}
\text{(a)}~~\lambda\left\|\xmin\right\|_2^2 \leq \lambda\left\|\xtrue\right\|_2^2+ \left\|A\xtrue-y\right\|_2^2-\left\|A\xmin-y\right\|_2^2.
\end{equation}
Let $\btrue=(\beta_1(\xtrue_1),\dots,\beta_n(\xtrue_n))=(|\xtrue_1|,\dots,|\xtrue_n|)$, and $\bmin=(\beta_1(\xmin_1),\dots,\beta_n(\xmin_n))$, as defined in \eqref{mcplasso_gen}. 
For any $z,\beta \in\R^n$, let $\|z\|_{1,\beta}=\sum_j \beta_i|z_i|$ be the $\beta$-weighted $\ell_1-norm$.
\begin{equation*}
\text{(b)}~~\lambda\left\|\xmin\right\|_{1,\bmin} \leq \lambda\left\|\xtrue\right\|_{1,\btrue}+ \frac{1}{2}\left\|A\xtrue-y\right\|_2^2-\frac{1}{2}\left\|A\xmin-y\right\|_2^2.
\end{equation*}
\end{lemma}
%
\begin{theorem}\label{robustheorem_gen}
Let $\xtrue\in\alf^n$ be a $k$-sparse signal with support $S$, and $\xtrue+\delta$ be its noisy version (\emph{i.e.}, the signal is not exactly sparse). Let $y=A(\xtrue+\delta)+\epsilon$, where $\epsilon$ is the measurement noise. Let us assume $d\geq 1$. If $A$ satisfies the RNSP of order $k$ with parameters  $\tau>0$ and $\rho\in(0,p^{-1})$ (where $p>0$ will be defined in the proof), then
$$\left\| \xmin - \xtrue\right\|_1 \leq C_1 \left\| A\delta+\epsilon\right\|_2^2+C_2\left(2\left\| A\delta+\epsilon\right\|_2+ qd\sqrt{\lambda  k}\right)$$
where $C_1$ and $C_2$ are positive constants assessed in the proof.
\end{theorem}
\begin{proof}
Let $h=\xmin - \xtrue$. 
\begin{equation*}
\left\| h_{S^c} \right\|_{1}=\left\| \xmin_{S^c} \right\|_1= \left\| \xmin \right\|_{1} -\left\| \xmin_S \right\|_{1}=\left\| \xmin \right\|_{1} -\left\| \xtrue_S+h_S \right\|_{1}.
\end{equation*}
For any $i\in S$, let $p_i\in\R$ such that $\bmin_i\xmin_i=\btrue_i\xtrue_i-p_i h_i$, and let $p=\max_i |p_i|$. We can then write:
\begin{equation*}
\begin{split}
&\left\| \xmin \right\|_{1} -\left\| \xtrue_S+h_S \right\|_{1}\leq \left\| \xmin \right\|_{1,\bmin} -\left\| \sum_{i\in S}\btrue_i\xtrue_i+p_i h_i \right\|_{1}\\
&\leq\left\| \xmin \right\|_{1,\bmin} -\left\| \xtrue_S\right\|_{1,\btrue}+p\left\|h_S \right\|_{1}.
\end{split}
\end{equation*}
Since $\left\| h_{S} \right\|_{1}\leq \sqrt{k}\left\| h_{S} \right\|_{2}$ \cite{matrix}, applying the RNSP we have $\left\| h_{S} \right\|_{1}\leq \rho \left\| h_{S^c} \right\|_{1}+\tau\sqrt{k} \left\| Ah \right\|_{2}$, under the assumption $\rho<p^{-1}$, we have:
\begin{equation*}
\begin{split}
&\left\| h_{S^c} \right\|_{1}\leq \left\| \xmin \right\|_{1,\bmin} -\left\| \xtrue\right\|_{1,\btrue} +p\rho \left\| h_{S^c} \right\|_{1}+p\tau\sqrt{k} \left\| Ah \right\|_{2}\\
&~~\Rightarrow \left\| h_{S^c} \right\|_{1}\leq \frac{1}{1-p \rho} \left(\left\| \xmin \right\|_{1,\bmin} -\left\| \xtrue\right\|_{1,\xtrue}+ p \tau\sqrt{k}\left\| Ah \right\|_{2}\right).
\end{split}
\end{equation*}
Therefore, 
\begin{equation*}
\begin{split}
\left\| h \right\|_1 &= \left\| h_S \right\|_1 +\left\| h_{S^c} \right\|_1 \leq (1+\rho)\left\| h_{S^c} \right\|_1+\tau\sqrt{k} \|Ah\|_2\\
&\leq \frac{1+\rho}{1-p \rho}\left[ \left\| \xmin \right\|_{1,\bmin} -\left\| \xtrue\right\|_{1,\btrue}\right]+\frac{p+1-p\rho}{1-p\rho}\tau\sqrt{k}\left\| Ah \right\|_{2}.
\end{split}
\end{equation*}
We now exploit Lemma \ref{doublereg_gen_noise}.(b):
\begin{equation}
\left\|\xmin\right\|_{1,\bmin} -\left\|\xtrue\right\|_{1,\btrue} \leq \frac{1}{\lambda}\left\|A\xtrue-y\right\|_2^2
\end{equation}
along with \eqref{eq:Ah} and $A\xtrue-y=A\delta+\epsilon$ to conclude that
{\small{$$ \left\| h \right\|_1\hspace{-1mm} \leq \hspace{-1mm} \frac{1+\rho}{\lambda(1-p \rho)}\left\| A\delta+\epsilon\right\|_2^2 +\frac{p+1-p\rho}{1-p\rho}\tau\sqrt{k}\left[2\left\| A\delta+\epsilon\right\|_2\hspace{-1mm}+\hspace{-1mm}\sqrt{\lambda}\left\| \xtrue\right\|_2\right]\hspace{-1mm}.$$}}

The thesis is then proven with $C_1=\frac{1+\rho}{\lambda(1-p \rho)}$ and $C_2=\frac{p+1-p\rho}{1-p\rho}\tau\sqrt{k}$, and applying $\left\| \xtrue\right\|_2\leq qd\sqrt{k}.$
\end{proof}
\begin{remark}\label{etalambda2}
Similarly to Remark \ref{etalambda}, if $\left\| A\delta+\epsilon\right\|_2\leq \eta $ for some known $\eta>0$, one can set $\lambda=\eta$, so that the error bound is of order $\sqrt{\eta}$ for small $\eta$, and, in particular, the error tends to zero when $\eta$ tends to zero.
\end{remark}

\section{MADMM: ADMM for MCP-Lasso }\label{sec:algorithm}
In this section, we present a novel algorithm, called MADMM, which is based on the Alternating Direction Method of Multipliers (ADMM, \cite{boy10}) and is designed for MCP-Lasso \eqref{mcplasso}. We derive it and we discuss its convergence. In the next section, numerical simulations demonstrate that MADMM outperforms the state-of-the-art methods.
\subsection{MADMM for $\fun$}
Following \cite[Section 2]{hon15}, we rewrite MCP-Lasso \eqref{mcplasso} as the following linearly constrained problem:
\begin{equation}\label{constrained}
\begin{split}
&\min_{x,z} \frac{1}{2}\left\|y-A x \right\|_2^2-\frac{\lambda}{2}\left\|x \right\|_2^2+\lambda d \left\|z\right\|_1\\
&\text{s.t. } z=x, x \in [-d,d]^n.
\end{split}
\end{equation}
We notice that the concave penalty has been split into a term depending on $x$ and another depending on $z$. This choice will allow us to prove the convergence. We know write the corresponding augmented Lagrangian (analogously to \cite[Section 2]{hon15}):
\begin{equation}\label{augmented_lagrangian}
\begin{split}
\lag(x,z)&=\frac{1}{2}\left\|y-A x \right\|_2^2-\frac{\lambda}{2}\left\|x \right\|_2^2+\lambda d \left\|z\right\|_1\\&~~+\mu^T(x-z)+\frac{\alpha}{2}\left\|x-z \right\|_2^2
\end{split}
\end{equation}
where $\mu$ is the dual variable, and $\alpha>0$. At this point, we can apply the classical ADMM procedure \cite[Section 2]{hon15}), which consists in the iteration (until convergence) of three steps: 1) minimization of $\lag$ with respect to $x$; 2) minimization of $\lag$ with respect to $z$; 3) update of the dual variable $\mu$.

The minima with respect to $x$ and $z$ can be easily computed in closed form. We notice that fixed $z$ (respectively, $x$) the problem is convex in $x$ (respectively, $z$). Since the solution must be in $[-d, d]^n$, given the convexity, it suffices to find the minima for $x$ and $z$ and then project them onto $[-d, d]^n$, that is, if $|x_{i}|>d$, then $x_{i}=\text{sign}(x_i)d$ (and analogously for $z$). We indicate by $P$ the projection onto $[-d, d]^n$.
The so-obtained procedure is written in Algorithm \ref{alg:MADMM}. $\soft_a:\R^n\to\R^n$, $a>0$, is the component-wise soft thresholding operator, defined as follows: $\soft_a(x)=0$ if $|x|<a$; $\soft_a(x)=x-a$ if $x>a$; $\soft_a(x)=x+a$ if $x<-a$.
\begin{algorithm}
\setstretch{1.3}
     \renewcommand{\algorithmicrequire}{\textbf{Input:}}
    \renewcommand{\algorithmicensure}{\textbf{Output:}}
  \caption{MADMM for $\fun$}\label{alg:MADMM}
  \begin{algorithmic}[1] 
   \REQUIRE $A,y, \lambda>0, \alpha>0$
 \ENSURE  $x_{T}$ = estimate of $\xtrue$
 \\ \textit{Notation}:$P$ = operator that projects onto $[-d, d]^n$; $\soft_a$ = soft thresholding operator
    \STATE Initialize $z_0=\mu_0=0\in\R^n$
    \FORALL{$t=1,\dots,T$}
    \STATE  $x_t=\argmin{x\in[-d,d]^n} \lag(x,z_{t-1})$\\
     $=P\left(\left[A^T A + (\alpha-\lambda)I\right]^{-1}\big( A^T y + \alpha z_{t-1} - \mu_{t-1}\big)\right)$
    %
    \STATE $z_t=\argmin{z\in[-d,d]^n} \lag(x_t,z)=P\left(\soft_{\frac{\lambda\beta}{\alpha}}\big(x_t+\frac{\mu_{t-1}}{\alpha}\big)\right)$
    %
    \STATE $\mu_t=\mu_{t-1}+\alpha(x_t-z_t)$
    \ENDFOR     
  \end{algorithmic}
\end{algorithm} 

\subsection{Convergence of MADMM for $\fun$}
While the convergence of ADMM for convex problems was established and largely studied some years ago \cite{boy10}, the convergence of ADMM for non-convex problems has been faced more recently \cite{wan15, li15, hon15, hon16}. For our purpose, we mainly rely on the convergence study in \cite{hon16}. 

In \cite{hon16}, the functional $\sum_{k=1}^K g_k(x)+h(x)$ has been considered, where each $k$ represents an agent assuming a distributed setting. In this work, the setting is centralized, \emph{i.e.}, $K=1$. Let us rewrite \cite[Assumption A]{hon16} in our centralized framework.
Consider the functional $g(x)+h(x)$, $x\in\X$, and the augmented Lagrangian:
$$\lag(x,z)=g(x)+h(z)+ \mu^T(x-z)+\frac{\alpha}{2}\left\|x-z \right\|_2^2.$$
\begin{assumption}\label{theassumption}

\begin{enumerate}
\item $\left\|\nabla g(x) - \nabla g(z)\right\|_2\leq C \left\|x-z \right\|_2$ for some $C>0$,  $\forall x,y\in\X$;
\item $h$ is convex (possibly non-smooth);
\item $\X$ is closed and convex;
\item the sub-problem $\min_x \lag(x,z)$ is strongly convex, with strongly convexity coefficient $\gamma$;
\item $\alpha \gamma> 2 C^2$ and $\alpha \geq C$;
\item $g(x) +h(x)$ is bounded from below.
\end{enumerate}
\end{assumption}
\begin{lemma}
MADMM for problem \eqref{constrained} satisfies Assumption \ref{theassumption}.
\end{lemma}
\begin{proof}
Let $g(x)=\frac{1}{2}\left\|y-A x \right\|_2^2-\frac{\lambda}{2}\left\|x \right\|_2^2$ and $h(z)=\lambda d\left\|z\right\|_1$.
\begin{enumerate}
	\item $\left\|\nabla g(x) - \nabla g(z)\right\|_2= \left\|(A^T A-\lambda I)(x-z)\right\|\leq C\left\|(x-z)\right\|_2$ where $C=\left\|A^T A-\lambda I\right\|_2$;
\item  $h$ is convex as $\ell_1$-norm is convex;
\item $[-d,d]^n$ is convex and closed;
\item The Hessian of 
$\frac{1}{2}\left\|y-A x \right\|_2^2-\lambda\frac{1}{2}\left\|x \right\|_2^2+ \mu^T(x-z)+\frac{\alpha}{2}\left\|x-z \right\|_2^2$ with respect to $x$ is $A^T A-\lambda I+\alpha I$, which is positive definite for any $\alpha>\lambda $;
\item Let $\gamma=\left\|A^T A-\lambda I+\alpha I\right\|_2$. For any sufficient large $\alpha$, it is easy to have $\alpha\gamma>2C$ and $\alpha>C$;
\item $\fun$ is non-negative.
\end{enumerate}
\end{proof}
\begin{theorem}\label{convergence_fun}
MADMM for $\fun$ converges to the set of stationary points of problem \eqref{constrained}, {\textcolor{black}{ \emph{i.e.}, if $x_t$ is the sequence generated by MADDM and $Z$ is the set of the stationary points, $\lim_{t\to\infty}\min_{z\in Z}\|x_t-z\|_2=0$.}}
\end{theorem}
\begin{proof}
Given Assumption \ref{theassumption}, \cite[Theorem 2.4- Point 3.]{hon16} guarantees that if $\X$ is compact, then the algorithm converges to the set of stationary points of \eqref{constrained}.
\end{proof}
%
\subsection{MADMM for $\gun$}
In Algorithm \ref{alg:MADMM_g}, we write the MADMM procedure for $\gun$. The unique substantial difference from Algorithm \ref{alg:MADMM} is that during each iteration also the feasibility set is updated, accordingly to the definition of $\gun$, in which the weights of the $\ell_1$-norm are $\beta_i(x_i)=\min\{\alpha\in\alf\text{ s.t. } |x_i|\leq \alpha\}$. In particular, the set size is non-increasing. 
\begin{algorithm}
 \renewcommand{\algorithmicrequire}{\textbf{Input:}}
    \renewcommand{\algorithmicensure}{\textbf{Output:}}
 \setstretch{1.3}
  \caption{MADMM for $\gun$}\label{alg:MADMM_g}
  \begin{algorithmic}[1] 
      \REQUIRE $A,y,\lambda>0, \alpha>0 $
 \ENSURE  $x_{T}$ = estimate of $\xtrue$
    \\ \textit{Notation}: $P_t$ = operator that projects onto $\mathcal{X}_t$; $\soft_a$ = soft thresholding operator
    \STATE Initialize $z_0=\mu_0=0\in\R^n, \beta_0=qd (1,\dots,1)^T\in\R^n$
    \FORALL{$t=1,\dots,T$}
    \STATE $\mathcal{X}_t=\prod_{i=1}^n [-\beta_{t-1,i}, \beta_{t-1,i}]$
    \STATE  $x_t=\argmin{x\in\mathcal{X}_t} \lag(x,z_{t-1})$\\
    $=P_t\left(\left[A^T A + (\alpha-\lambda)I\right]^{-1}\big( A^T y + \alpha z_{t-1} - \mu_{t-1}\big)\right)$
    %
    \STATE $z_t=\argmin{z\in\mathcal{X}_t} \lag(x_t,z)=P_t\left(\soft_{\frac{\lambda\beta_t}{\alpha}}\big(x_t+\frac{\mu_{t-1}}{\alpha}\big)\right)$ 
    %
    \STATE $\beta_{t,i}= \sum_{j=1}^q dj \I_{(d(j-1), dj]}(|z_{t,i}|)$, $i=1,\dots, n$

    \STATE $\mu_t=\mu_{t-1}+\alpha(x_t-z_t)$
    \ENDFOR     
  \end{algorithmic}
\end{algorithm} 

\subsection{Convergence of MADMM for $\gun$}
\begin{theorem}\label{convergence_gun}
MADMM for $\gun$ converges to the set of stationary points.
\end{theorem}
\begin{proof}
For any $i\in\{1,\dots, n\}$, the discrete sequence $\{\beta_{t,i}\}_{t=1,2,\dots}$ is non-increasing and lower bounded by zero, hence it converges in a finite number of iterations. After the stabilization of $\beta$, we are in the same setting of MADMM for $\fun$, hence convergence can be proven as in Theorem \ref{convergence_fun}.
\end{proof}

\subsection{How to check if the solution is exact}
In the previous section, we have shown that, under mild conditions, the global minimum of MCP-Lasso is the original signal in the noise-free case. This never happens with Lasso, where the minimum is always affected by a bias proportional to $\lambda$. This theoretical advantage of MCP-Lasso, however, is at the price of non-convexity: the minimum of Lasso can be achieved straightforwardly leveraging on  convexity, which is not true for MCP-Lasso. The problem is then how to achieve the global minimum of MCP-Lasso. 

An important help in this direction comes from the following proposition (valid for the ternary alphabet; the extension to the generic alphabet is left for future work).
\begin{proposition}\label{howtocheck}
If ker$(A_S^T A)\notin \Z$, $\xtrue$ is the unique point in $\{0,\pm d\}^n$ that can be a stationary point of MADMM.
\end{proposition}
\begin{proof}
It is easy to check that a stationary point $x$ of MADMM satisfies the property $A^T(Ax-y)=\lambda x -\mu$,
where $\mu\in [-\lambda,\lambda]^n$, and more precisely $\mu_i=\lambda \text{sign}(x_i)$ when $x_i\neq 0$, and $\mu_i\in(-\lambda,\lambda)$ when $x_i=0$.

Let $x\in \{0,\pm\beta\}^n$, $x\neq \xtrue$. For any $x_i\neq 0$, we must have $A_S^T A (x-\xtrue)=0$, where $x-\xtrue\in d\{0,\pm 1,\pm 2\}$. However this contradicts the kernel hypothesis.

Otherwise, if $x=\xtrue$,  clearly $A_S^T A (x-\xtrue)=0$, and $\mu_i=0$ when $x_i=0$.
\end{proof}

The hypothesis on the kernel is quite similar to the general position property \cite{tib13} and is almost always satisfied in practical situations. For instance, it is satisfied with probability 1 when the entries drawn from a continuous probability distribution \cite[Lemma 4]{tib13}. In the practice, we can exploit Proposition \ref{howtocheck} as follows: if MADMM finds a solution in $\{0,\pm d\}^n$, then this solution is exact. Moreover, our numerical experiments will show that MADMM is very fast, therefore, when a non-exact solution is found, it is feasible to re-run the algorithm with different initialization to look for the right one, as we will show in Section \ref{sec:numerical}. {\textcolor{black}{We remark that this does not guarantee to achieve the global minimum, therefore our approach is still sub-optimal. However, the experiments in Section \ref{sec:numerical} will show that re-running the algorithm with random initialization in many cases is sufficient to find the global minimum. A more rigorous search of the global minimum is beyond the scope of this work. We however mention that, since the proposed cost functionals are semi-algebraic, the approach proposed in \cite{las15} might be considered: this methodology guarantees the global optimization, at the price of a high computational complexity. A thorough comparison between our sub-optimal approach and semi-algerbraic solutions will be proposed as future work.}}

As already said, an other algorithm to solve Lasso is IST \cite{dau04}, whose convergence is easy to prove \cite{for10} if compared to ADMM. Given its simplicity, IST is widely used, in particular in the distributed context \cite{rfm15, fox16}. However, it is well known that, in the practice, ADMM converges in a significantly lower number of iterations than IST. With this  motivation, in this paper we do not consider  IST, even though an IST formulation for MCP-Lasso can be easily derived. Moreover, in the practice, IST for MCP-Lasso converges, but the theoretical proof is not straightforward. We remind that \cite{woo16, bayram16} prove the convergence of IST in non-convex sparse problems, but limited to the case of positive definite $A^T A$ (which is not the CS case).

\section{Numerical experiments}\label{sec:numerical}
In this section, we show the efficiency of MADMM in terms of recovery accuracy, number of measurements, and speed of convergence through numerical simulations\footnote{{\textcolor{black}{The code to reproduce these simulations is available at https://github.com/sophie27/Recovery-of-sparse-finite-valued-signals}}}. First, we consider synthetic random signals; second, we tackle a localization application.
\subsection{Random finite-valued signals}
In the first set of experiments, we consider $k$-sparse signals in $\alf^n$, where $n=100$ and $\alf=\{0,\pm 1\}$ or $\alf=\{0,\pm 1, \dots, \pm 5\}$. Support and non-zero values are chosen uniformly at random. We consider Gaussian sensing matrices $\mathcal{N}\sim (0,\frac{1}{m})$. We implement two versions of the proposed algorithm MADMM: the original one, which stops when a stationary point is achieved, and the "reshuffling" one, which leverages Proposition \ref{howtocheck} to check whether the achieved point is the desired signal. If not, the algorithm reshuffles the initialization point and repeats the procedure to search a solution in $\alf^n$ (or at least a solution closer to $\alf^n$). In these experiments, we reshuffle by reinitializing  $z$ uniformly at random in conv$(\alf)^n$. The algorithm is definitely stopped when the relative square distance from the found solution and $\alf^n$ is smaller than $10^{-4}$. We indicate this second version as MADMM-R. {\textcolor{black}{We set $\lambda=10^{-2}$, this value being chosen because, with high probability, it is smaller than the lowest positive eigenvalue of $A^TA$ (see Remarks \ref{qf}, \ref{qf2}, and Corollary \ref{cor:noncs}) when $n=100$ and $m<n$. For simplicity, we do not re-compute $\lambda$ for each $A$.}}

Given the estimate $\widehat{x}$, we consider two performance metrics: the relative square error (RSE) defined as $\|\xtrue-\widehat{x}\|^2/\|\xtrue\|^2$ and the count of exact recovery occurrences, that is, the number of experiments where $\widehat{x}= \xtrue$. In the first experiments we compare our approach to Lasso \cite{kei17} and BP \cite{fli18}, both solved via ADMM (we recall that BP can be used only in the noise-free setting). As Lasso solution has a bias, we finally project it onto $\alf^n$. 

We fix the ADMM parameter $\alpha=1$ for all the algorithms, and we stop them when the distance between two successive estimates is below $10^{-12}$. The results are averaged over 500 runs.
\begin{figure*}[ht]
\centering
\includegraphics[width=0.31\textwidth]{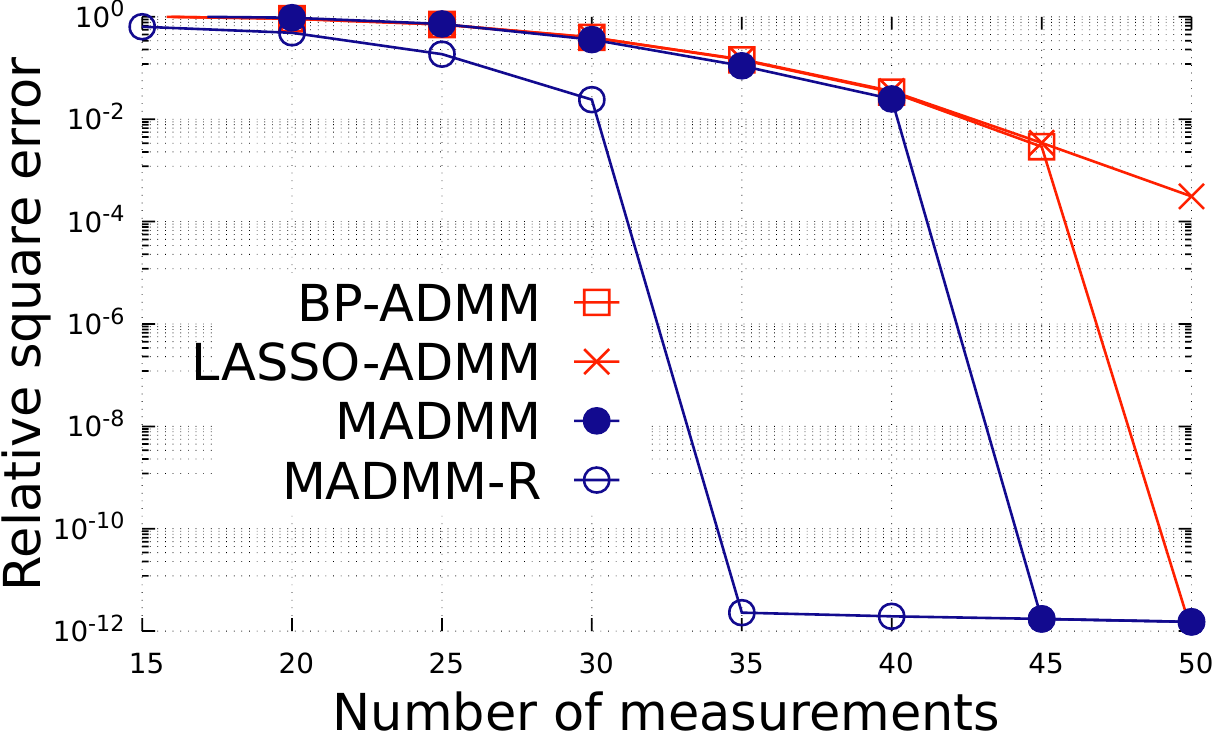}\,
\includegraphics[width=0.31\textwidth]{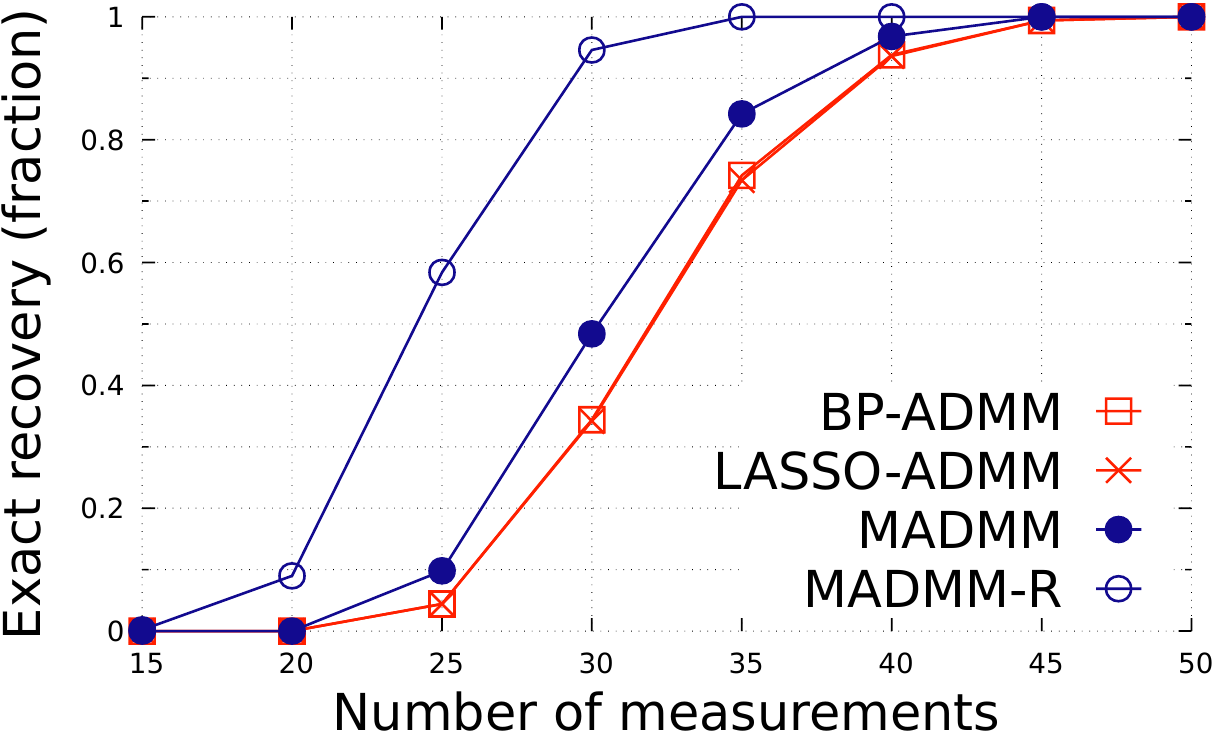}\,
\includegraphics[width=0.31\textwidth]{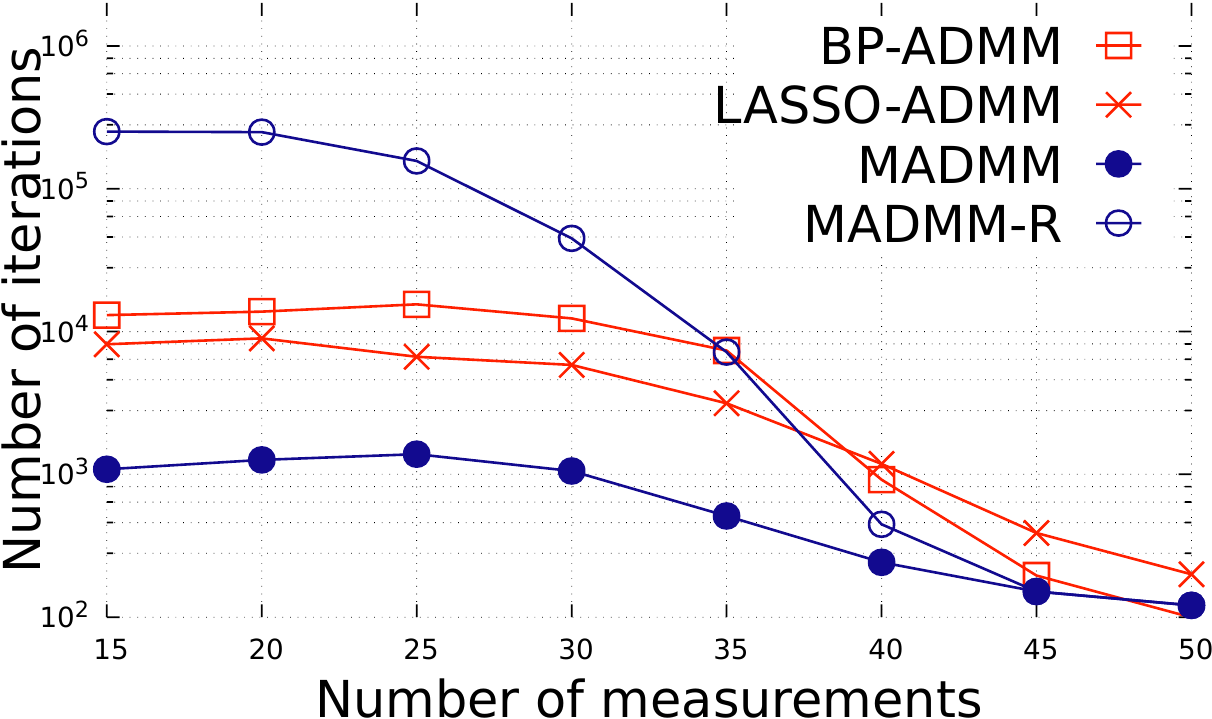}
\caption{$\alf=\{0,\pm 1\}$, $n=100$, $k=10$, $\lambda=10^{-2}$, noise-free case, mean over 500 runs.}\label{fig:1m}
\end{figure*}
\begin{figure*}[ht]
\centering
\includegraphics[width=0.31\textwidth]{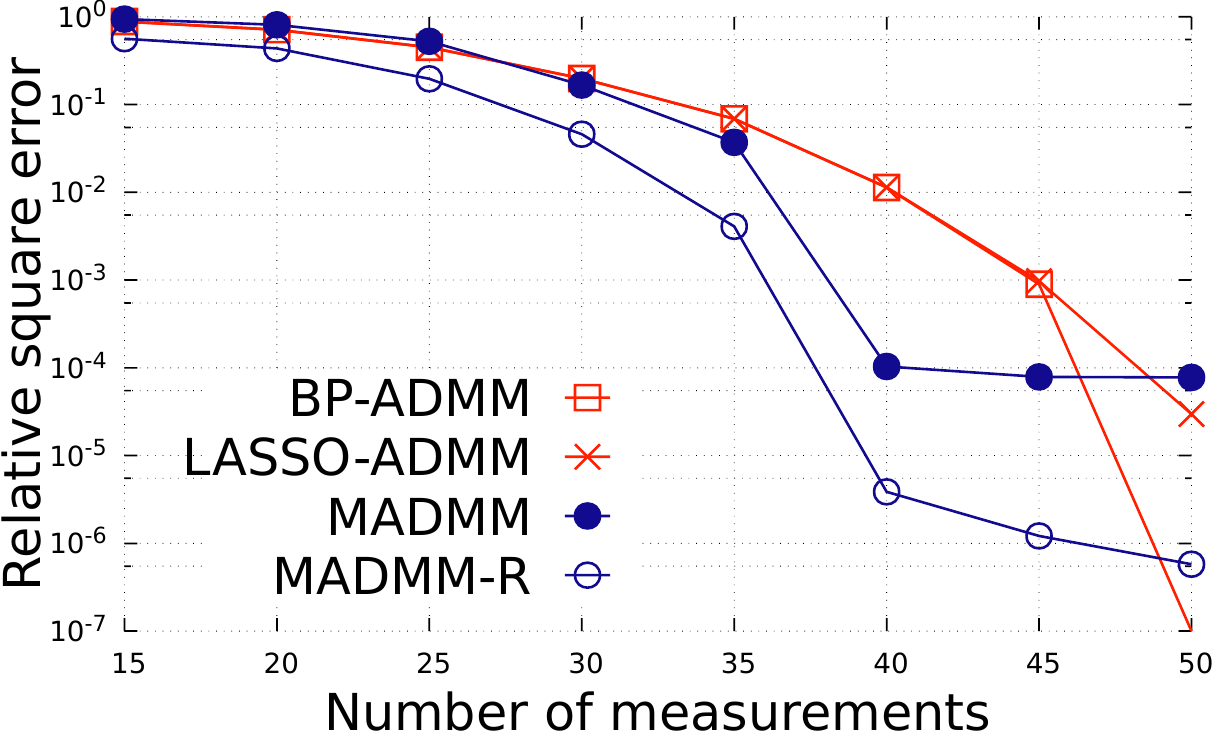}\,
\includegraphics[width=0.31\textwidth]{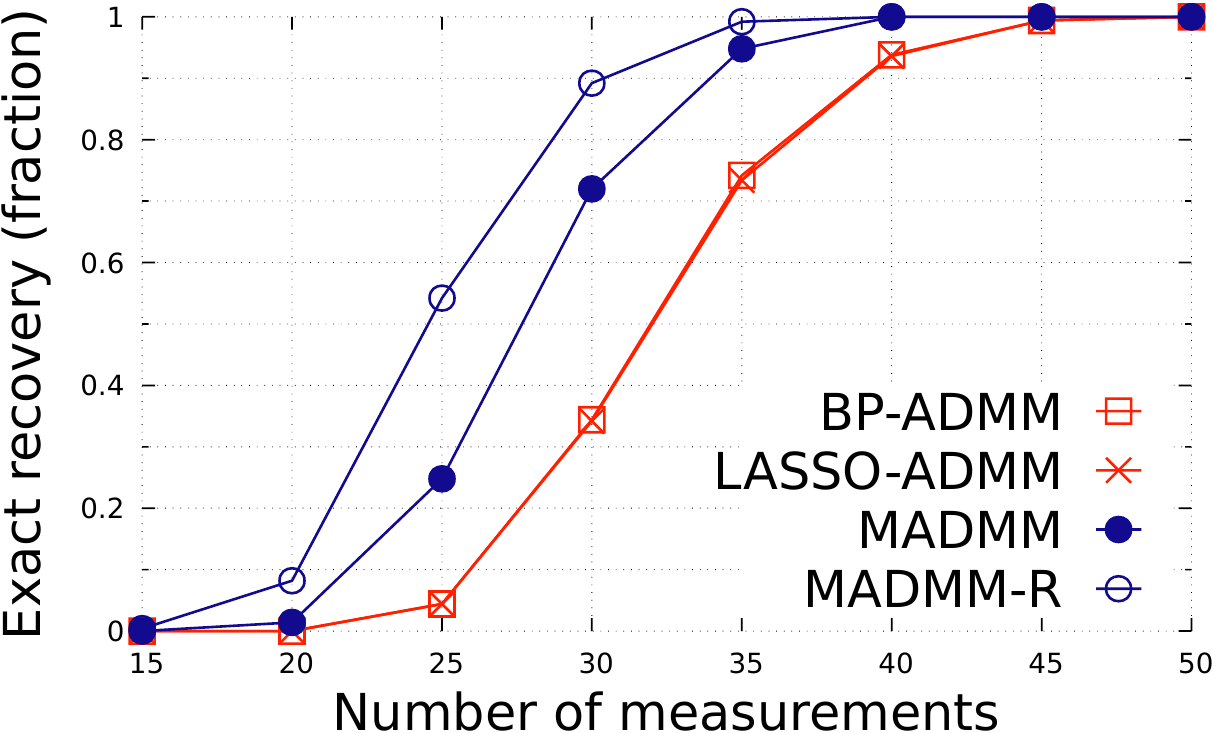}\,
\includegraphics[width=0.31\textwidth]{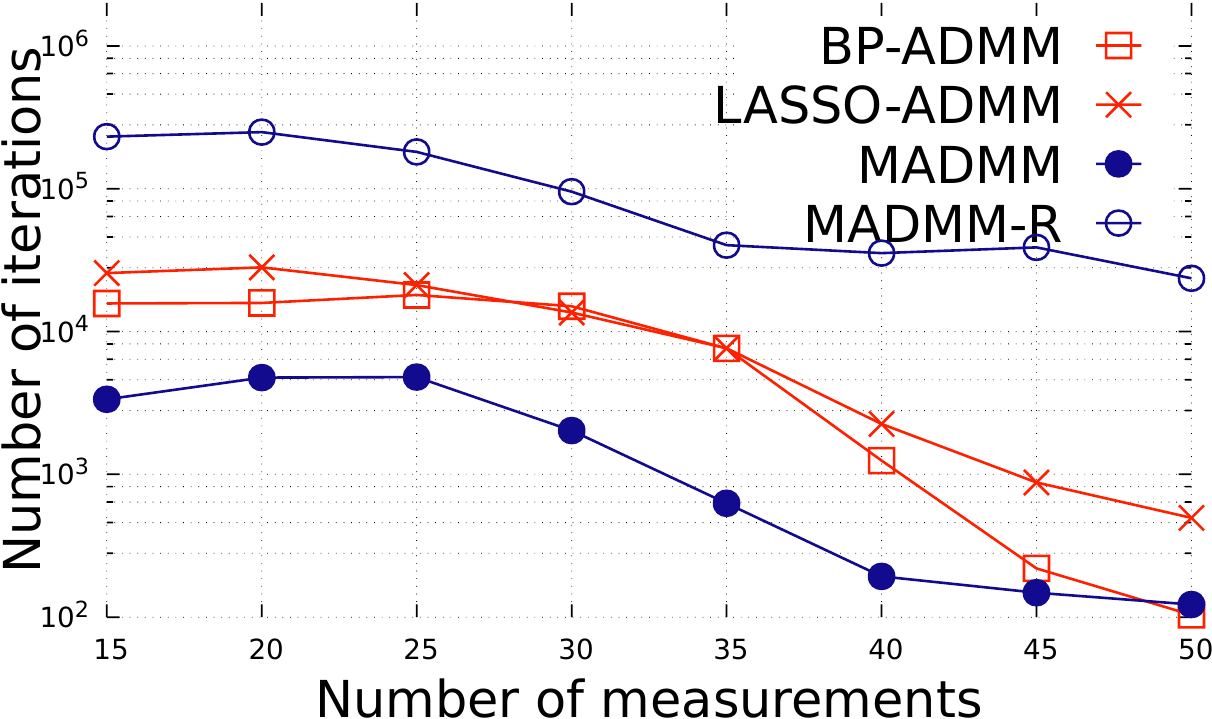}
\caption{$\alf=\{0,\pm 1, \dots, \pm 5\}$, $n=100$, $k=10$, $\lambda=10^{-2}$, noise-free case, mean over 500 runs.}\label{fig:5m}
\end{figure*}

In Figure \ref{fig:1m}, we show the results (accuracy metrics and number of iterations) of the experiments with $\alf=\{0,\pm 1\}$ in the noise-free setting; the sparsity is fixed to $k=10$ and we vary the number of measurements $m$. We remark that for each iteration step the complexity of ADMM and MADMM is actually the same, than the number of iterations determines the convergence time. We can appreciate that MADMM is more accurate (in both metrics) and quicker than ADMM. The accuracy can be further improved with MADMM-R (which always obtain the exact recovery with only $m=35$ measurements), at the price of a larger convergence time for smaller $m$'s.

A similar behavior can be observed for  $\alf=\{0,\pm 1, \dots, \pm 5\}$ (Figure \ref{fig:5m}). In this case, the improvement of MADMM with respect to ADMM is higher. MADMM-R instead has always has a larger convergence time, which means that a lot of reshuffling might be necessary to decrease RSE for larger alphabets. We recall however that the stopping criterion for reshuffling is set on the non-quantized final estimate; performing quantization we could expect an exact solution after less reshuffling iterations.

\begin{figure*}[ht]
\centering
\includegraphics[width=0.31\textwidth]{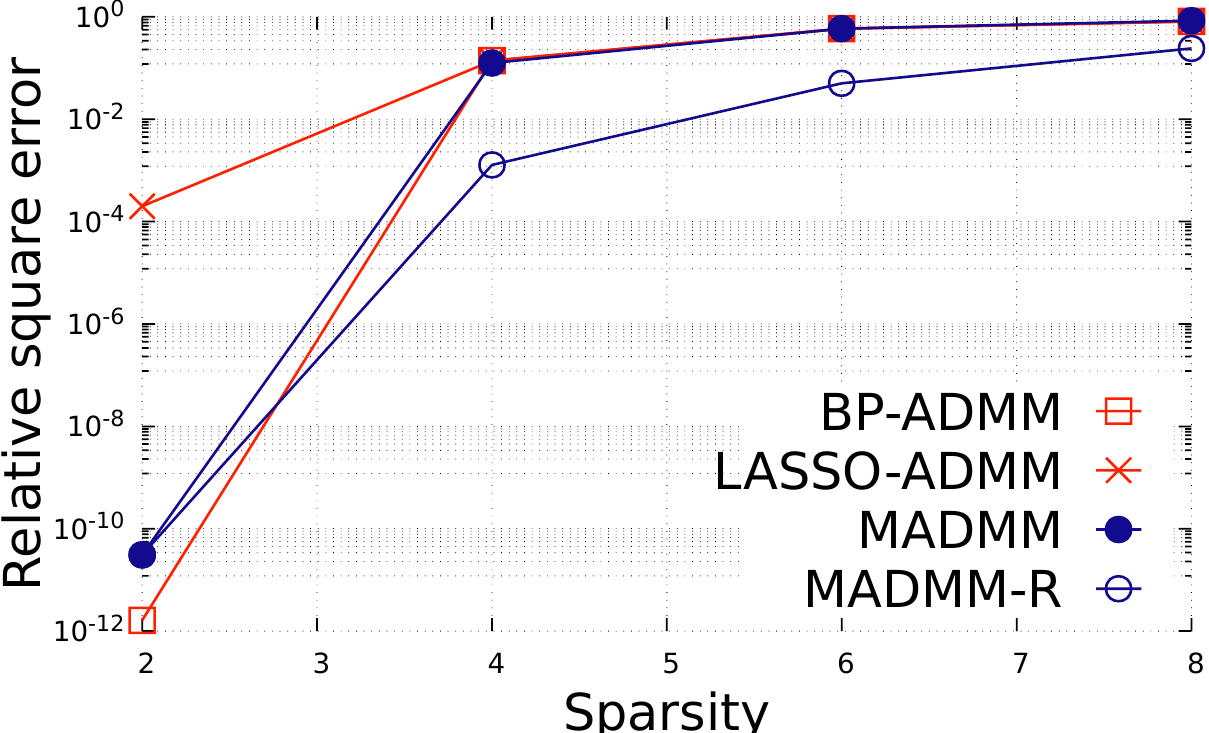}\,
\includegraphics[width=0.31\textwidth]{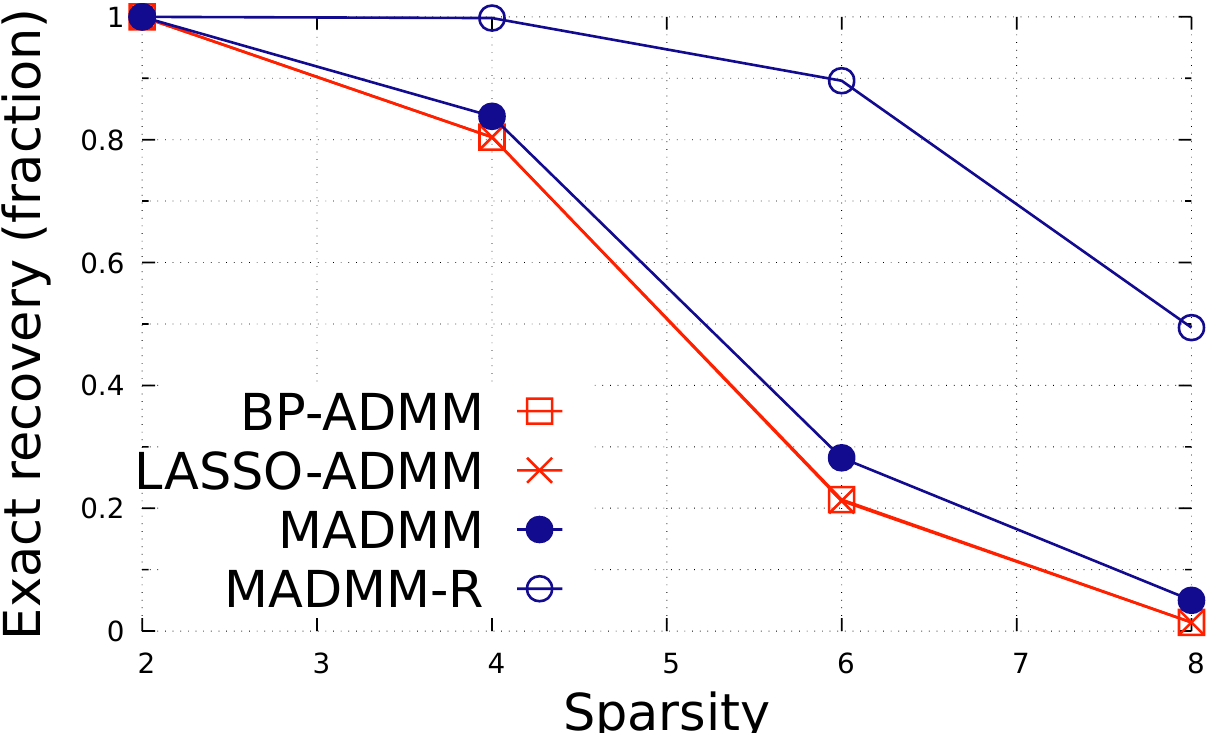}\,
\includegraphics[width=0.31\textwidth]{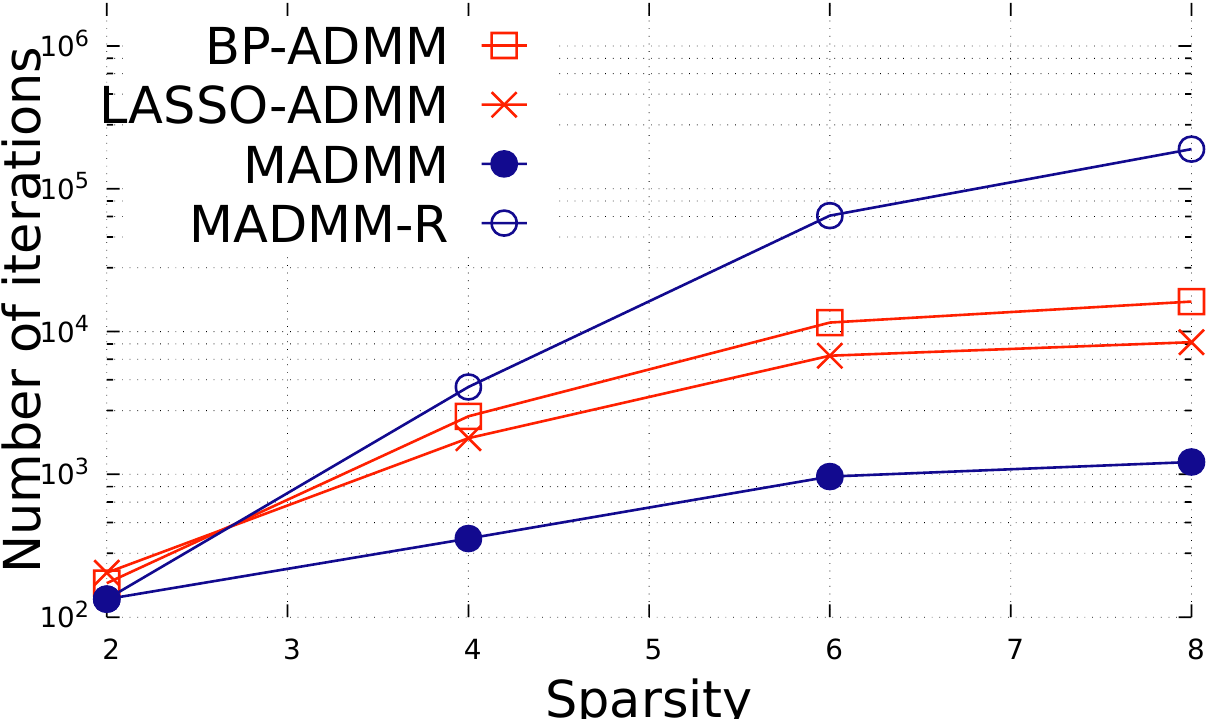}
\caption{$\alf=\{0,\pm 1\}$, $n=100$, $m=20$, $\lambda=10^{-2}$, noise-free case, mean over 500 runs.}\label{fig:1k}
\end{figure*}
\begin{figure*}[ht]
\centering
\includegraphics[width=0.31\textwidth]{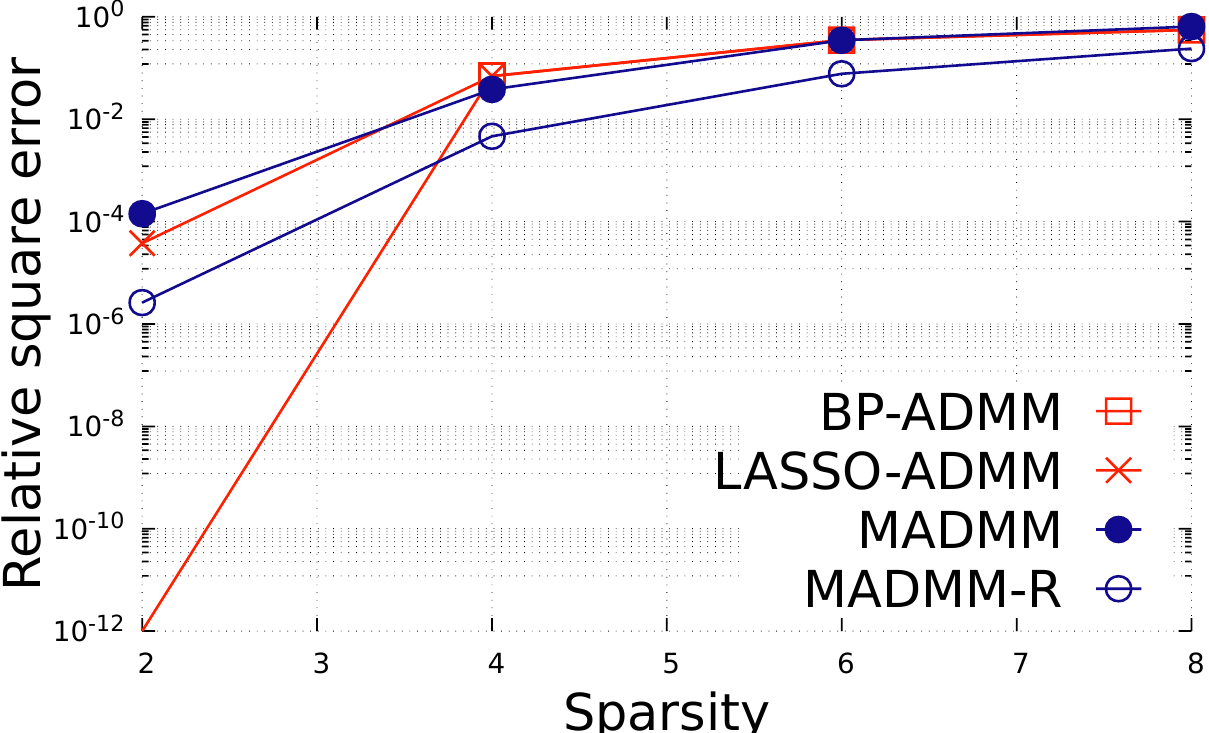}\,
\includegraphics[width=0.31\textwidth]{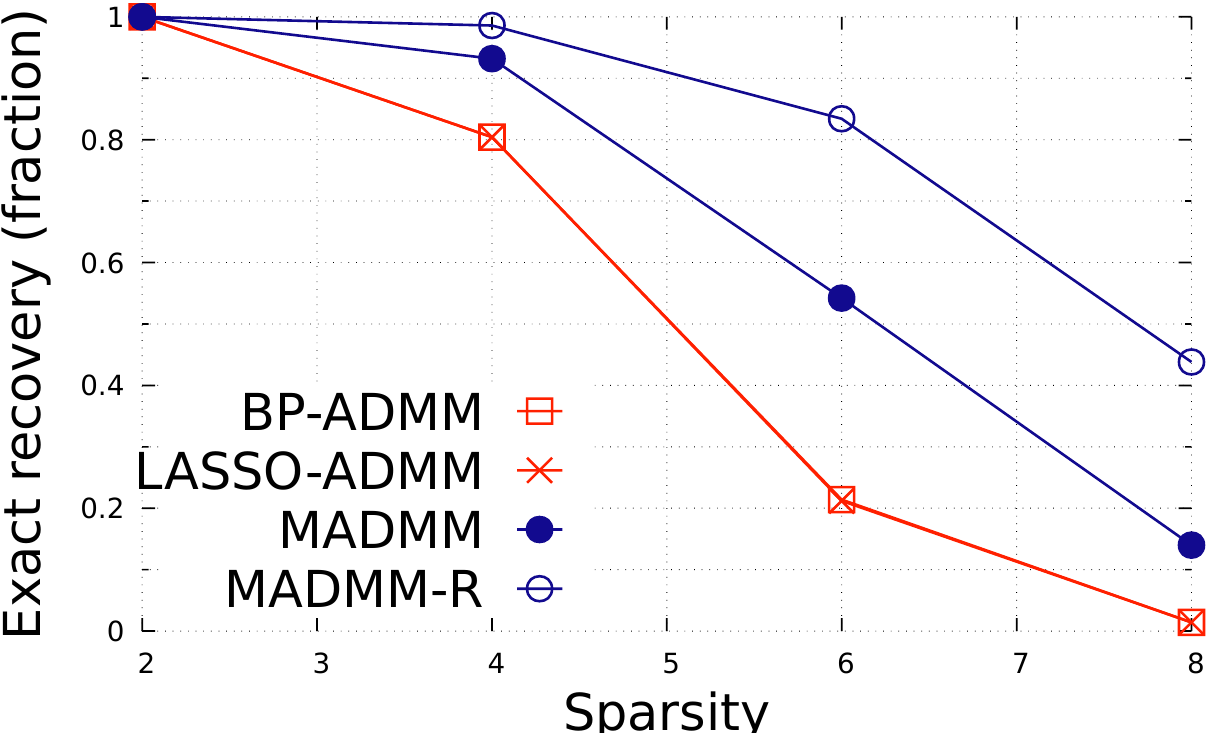}\,
\includegraphics[width=0.31\textwidth]{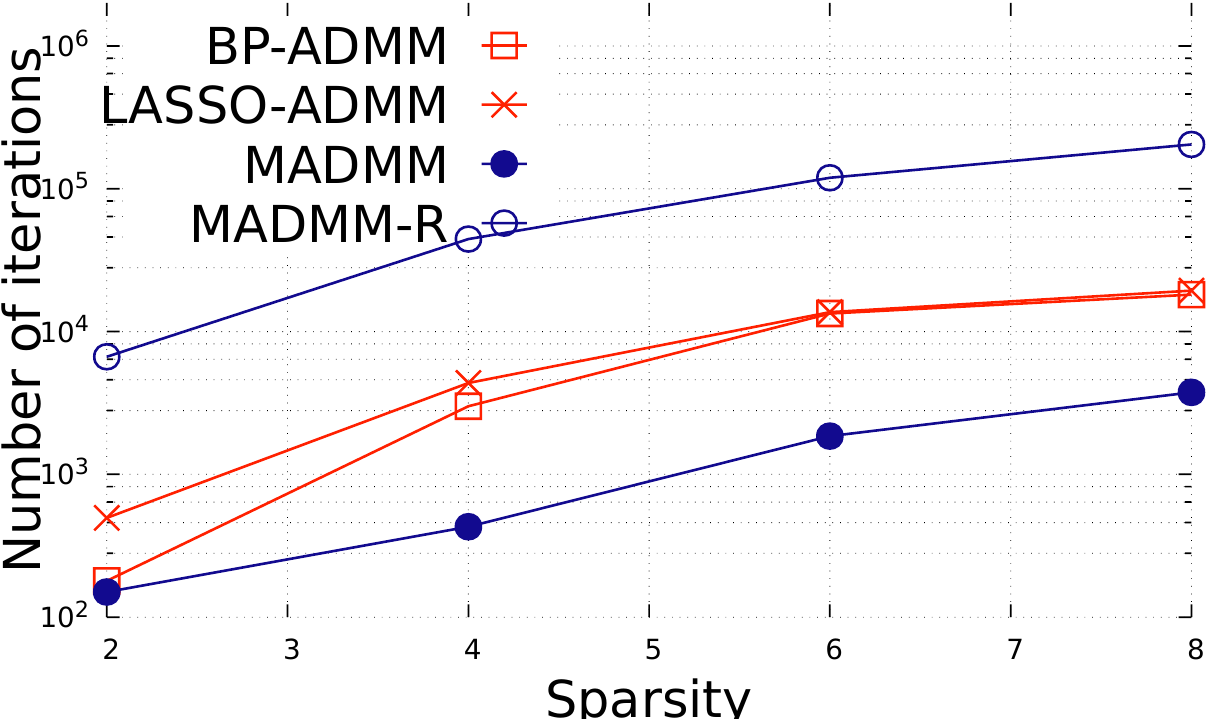}
\caption{$\alf=\{0,\pm 1, \dots, \pm 5\}$, $n=100$, $m=20$, $\lambda=10^{-2}$, noise-free case, mean over 500 runs.}\label{fig:5k}
\end{figure*}
{\textcolor{black}{In Figures \ref{fig:1k} and \ref{fig:5k}, we show the performance when the sparsity level $k$ varies, while $m=20$. The observed behaviors are in line with the previous results: MADMM and MADMM-R are more accurate than Lasso and BP, and MADMM is the fastest choice.}}

\begin{figure*}[ht]
\includegraphics[width=0.31\textwidth]{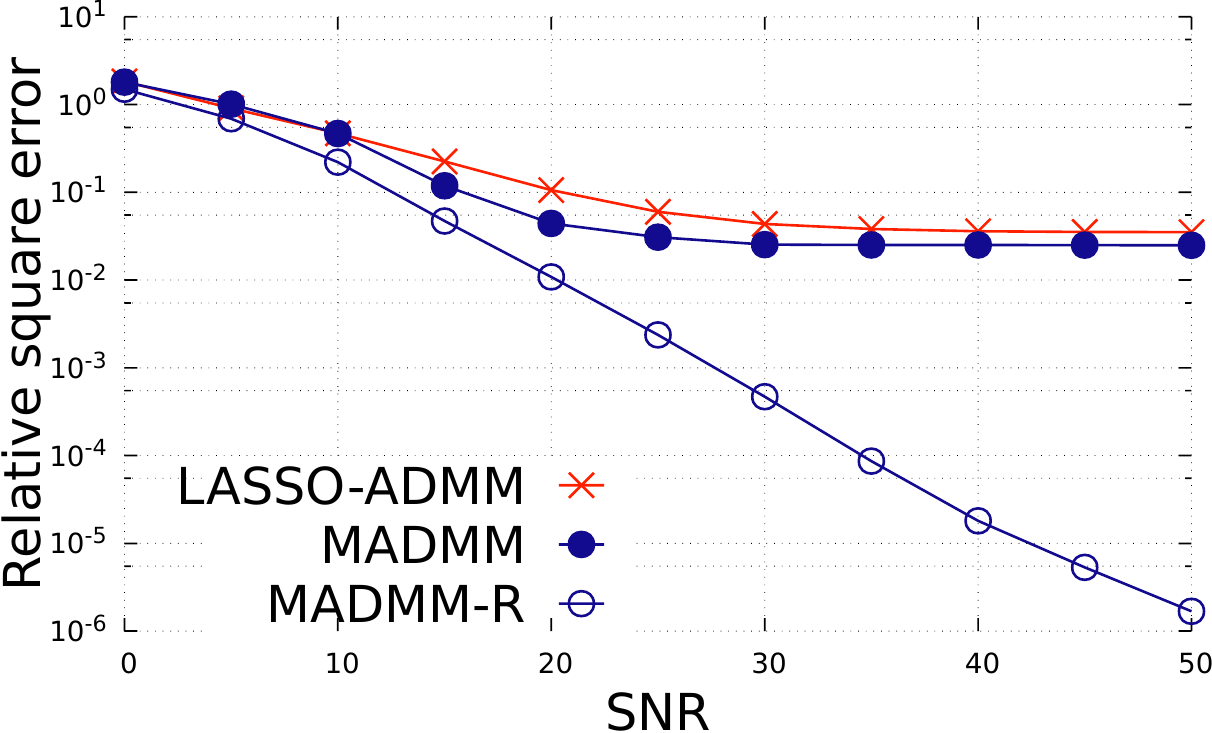}\,
\includegraphics[width=0.31\textwidth]{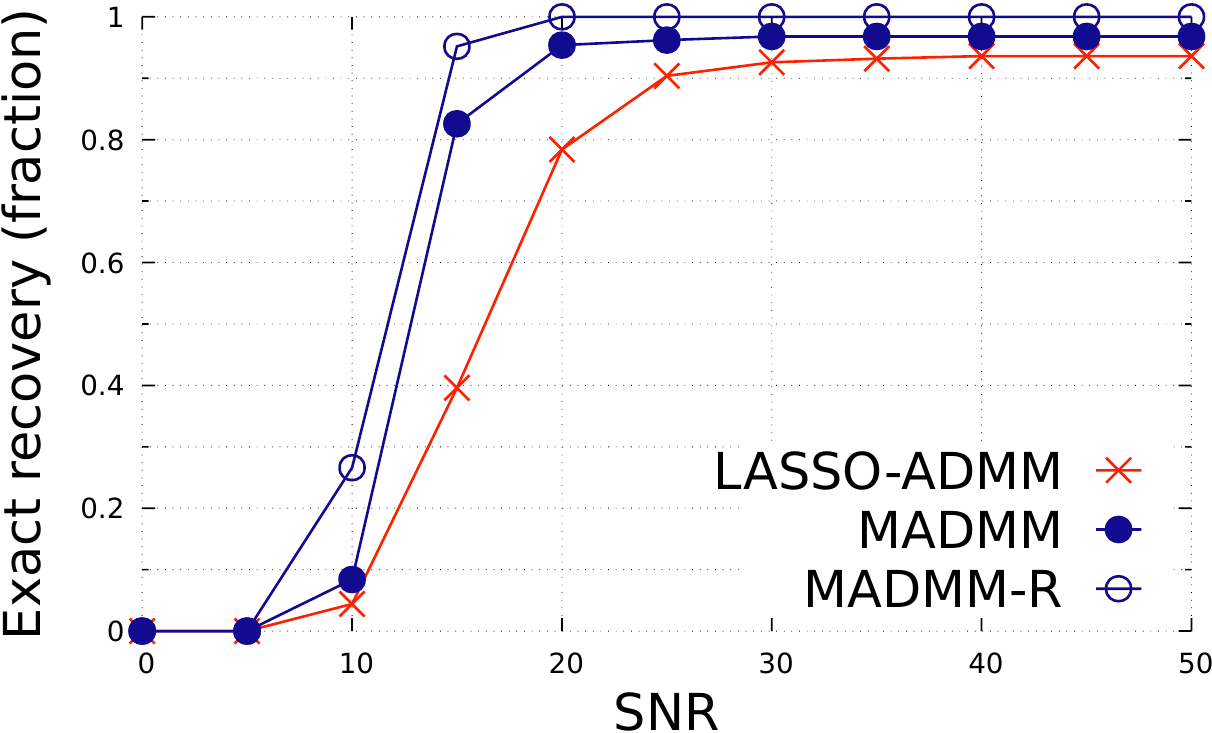}\,
\includegraphics[width=0.31\textwidth]{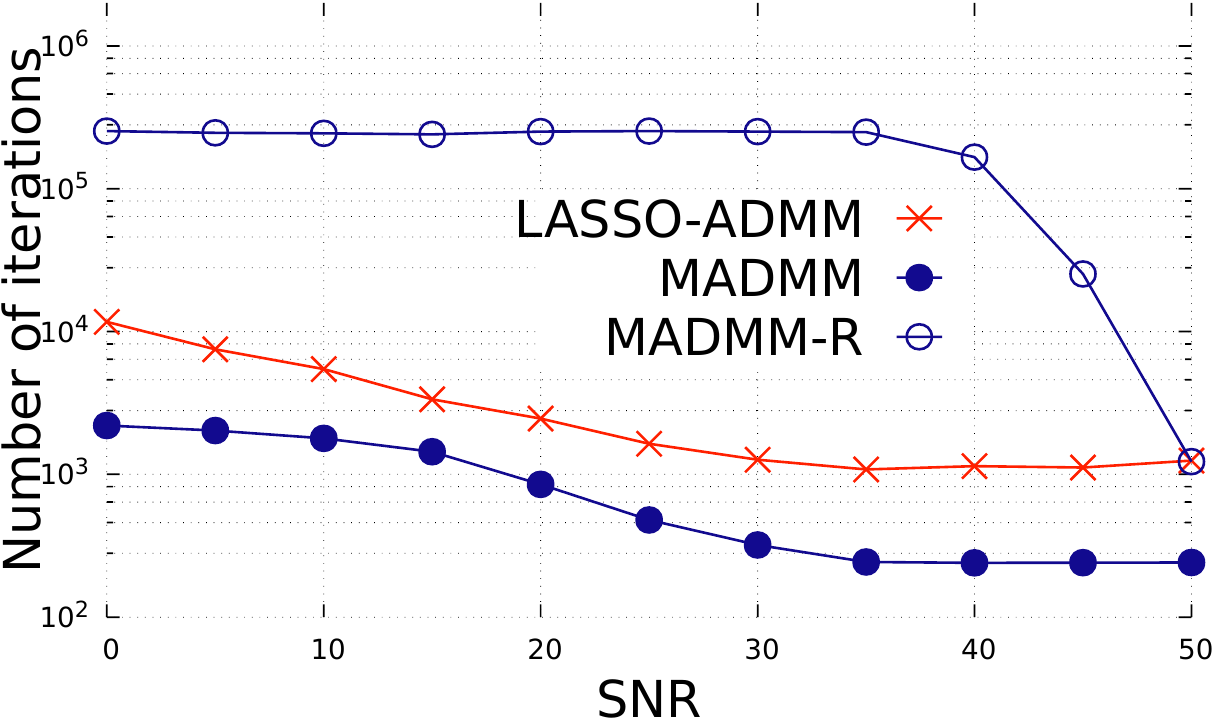}
\caption{$\alf=\{0,\pm 1\}$, $n=100$, $k=10$, $m=40$, $\lambda=10^{-2}$, noisy measurements, mean over 500 runs. SNR is expressed in dB.}\label{fig:1snr}
\end{figure*}

\begin{figure*}[ht]
\includegraphics[width=0.31\textwidth]{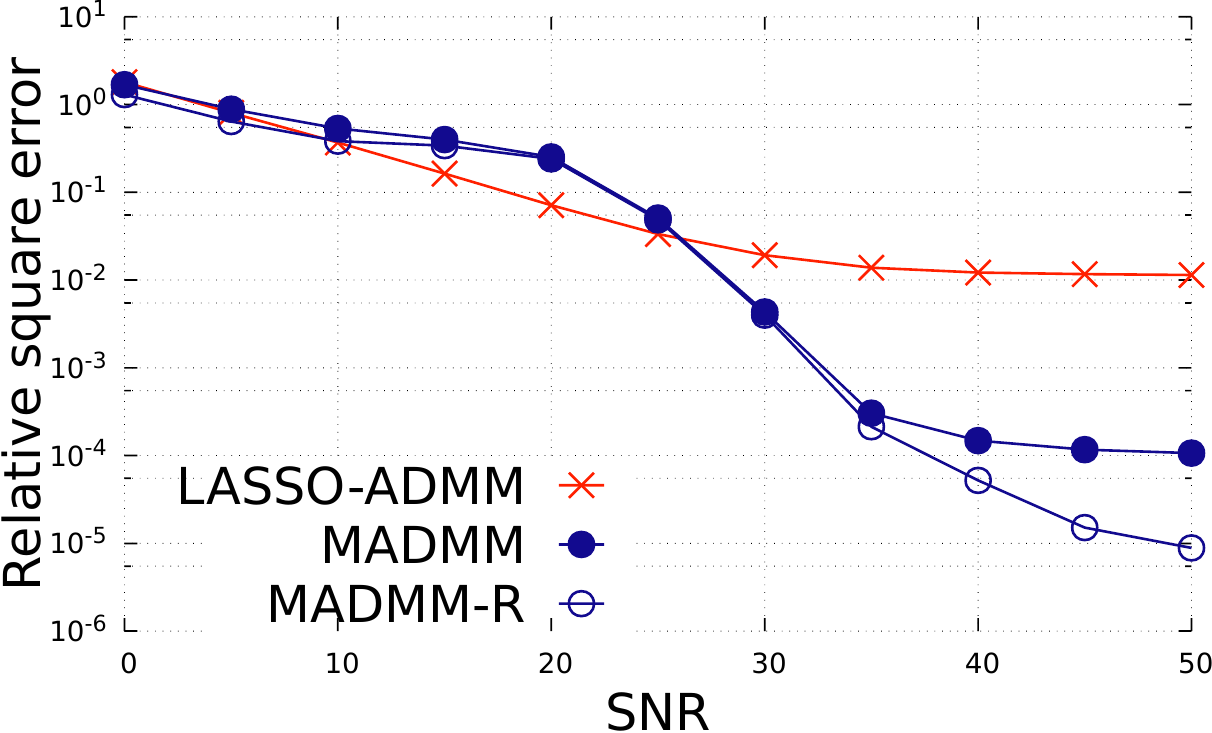}\,
\includegraphics[width=0.31\textwidth]{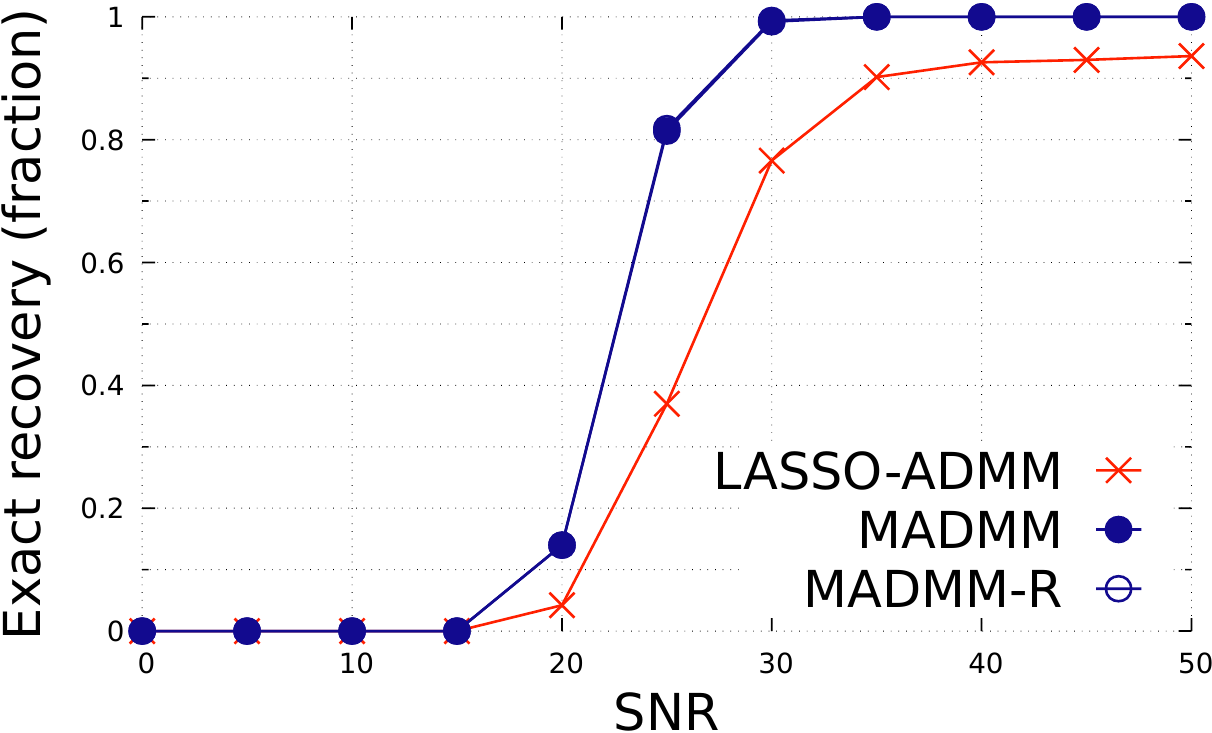}\,
\includegraphics[width=0.31\textwidth]{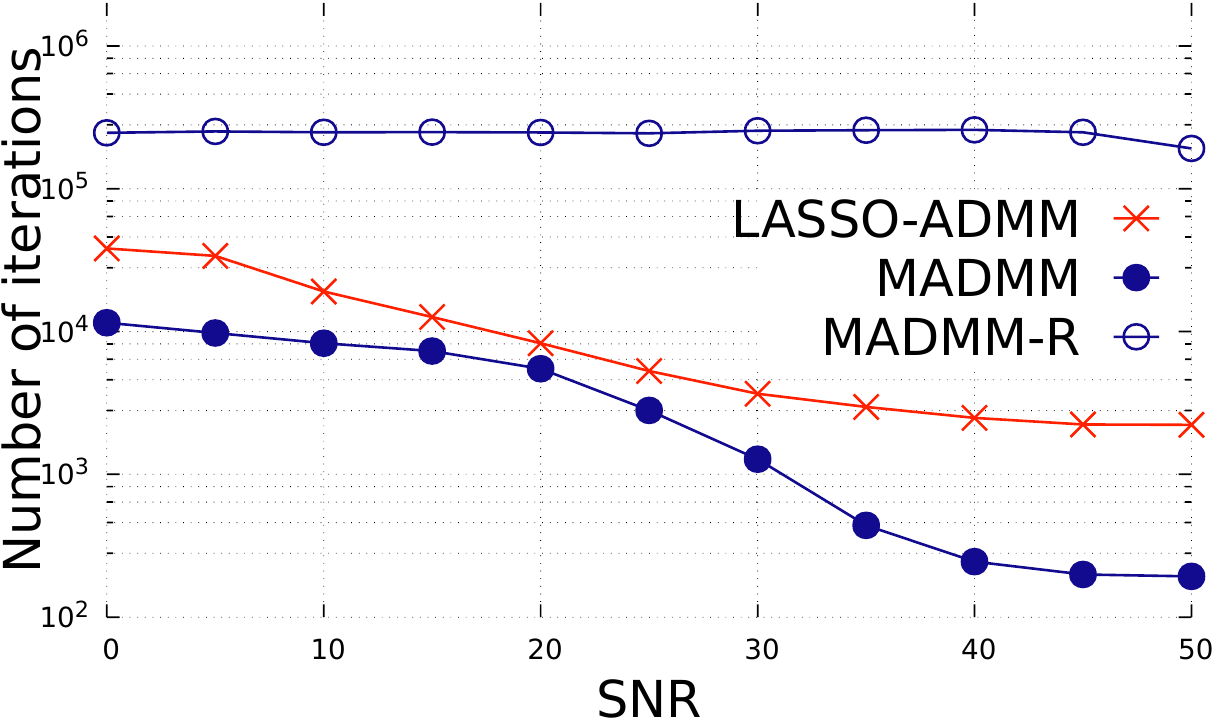}
\caption{$\alf=\{0,\pm 1, \dots, \pm 5\}$, $n=100$, $k=10$, $m=40$, $\lambda=10^{-2}$, noisy measurements, mean over 500 runs. SNR is expressed in dB.}\label{fig:5snr}
\end{figure*}

In the second set of experiments, we fix $k=10$, $m=40$, and we add some measurement Gaussian noise. In Figures \ref{fig:1snr} and \ref{fig:5snr}, we show the performance for different signal-to-noise ratios (SNR). Again, MADMM is observed to be more accurate and quicker than ADMM. For $\alf=\{0,\pm 1\}$, an evident gap is obtained at SNR$=15$dB, where ADMM recovers exactly only in $40\%$ of runs, while MADMM overcomes  $80\%$. Moreover, with reshuffling, we always obtain the exact solution for SNR$\geq 20$dB. For $\alf=\{0,\pm 1, \dots, \pm 5\}$ we observe that MADMM and MADMM-R have the same exact recovery rate (the higher number of iterations of MADMM-R is just to reduce the RSE). 

We conclude that MADMM generally performs better than Lasso and BP (solved via ADMM), in terms of recovery accuracy and computation complexity. Moreover, MADMM-R can be used when higher precision is required, at the price of slower convergence.

\begin{figure}[ht]
\centering
\includegraphics[width=0.48\textwidth]{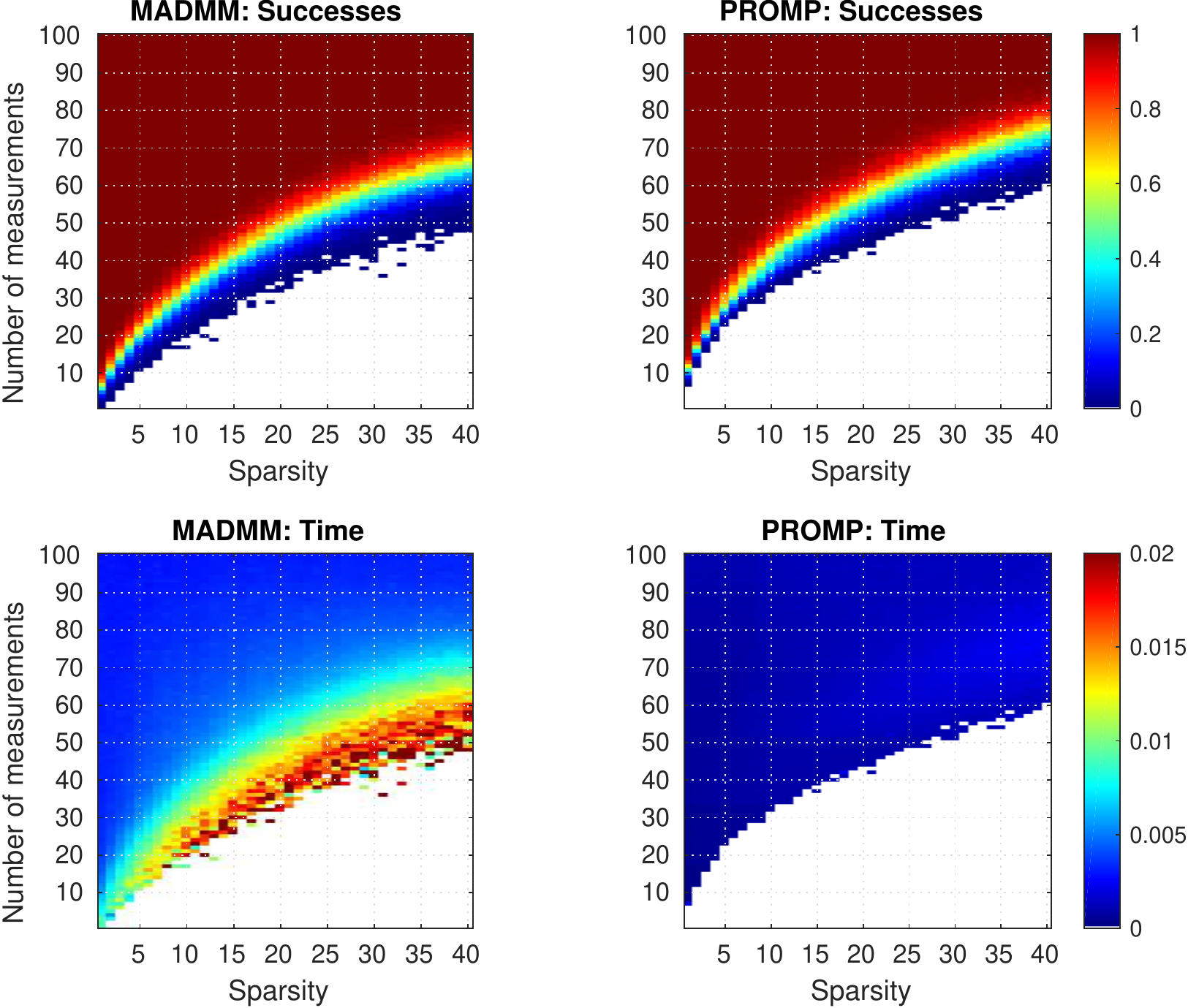}
\caption{MADMM vs PROMP \cite{fli18}, $n=100$, $\alf=\{0,\pm 1\}$, noise-free case, 500 runs. The run time is expressed in seconds.}\label{fig:promp}
\end{figure}

{\textcolor{black}{We finally compare MADMM with the latest algorithm PROMP \cite{fli18}, designed for lattice-valued signals. For this comparison, we consider a simulation setting proposed in \cite{fli18}: $A$ is Gaussian, $\alf=\{0,\pm 1\}$, $n=100$; no noise is added. In Figure \ref{fig:promp}, we show the rate of exact estimations over 500 trials, and the mean run time (in seconds). We can appreciate that MADMM achieves better recovery performance: less measurements are required to have success; in particular, the threshold to obtain $100\%$ of successes is around 10 measurements lower for MADMM. Concerning the run time, PROMP time is always of order $10^{-3}$ seconds, while MADMM generally requires less than $10^{-2}$ seconds in the $100\%$ successes zone, and is of order $10^{-2}$ in the other zones. This behavior is in line with what is observed in classical (non discrete-valued) CS: greedy methods are faster, while $\ell_1$-based methods require less measurements}}

\subsection{Multiple target localization via CS}
We now show the efficiency of MADMM in a practical problem, namely a multiple target localization problem \cite{fen09}. 
We consider a $20\times 20~\text{m}^2$ area, subdivided into $n=100$ cells of dimension $2\times 2~\text{m}^2$, and we simulate the following setting. $m<n$ sensors are deployed uniformly at random over the area. In the training phase, a target is placed in turn in each cell, and the corresponding received signal strength (RSS) at each sensor is measured, according to the IEEE 802.15.4 standard \cite[Equation 11]{fen09}, with SNR$=25$dB. Each sensor takes only one measurement, then the number of measurements is equal to the number of sensors (more measurements for each sensor  could be considered to improve the localization). In this way, we build the dictionary $A$.

Given $A$, the localization problem can be interpreted as the recovery of a binary signal $x\in\{0,1\}^n$ from $y=Ax$. Specifically, $x_i=1$ when a target transmits from cell $i$, and $x_i=0$ when the cell $i$ is empty. Even in the case of multiple targets, the number of targets is generally much smaller than the number of cells, which guarantees sparsity conditions. The matrix $A$ is deterministic and may not have sufficient properties of incoherence to apply CS. However, in \cite[Proposition 1]{fen09} it has been proven that after feasible orthogonalization, the problem becomes suitable for $\ell_1$-minimization. 

For this experiment we compare MADMM and Lasso solved via ADMM. We set $\lambda=10^{-3}$ and we stop the algorithms when the distance with previous step estimate is lower than $10^{-8}$. For MADMM, we clearly set $d=1$. 
\begin{figure}
\centering
\includegraphics[width=0.24\textwidth]{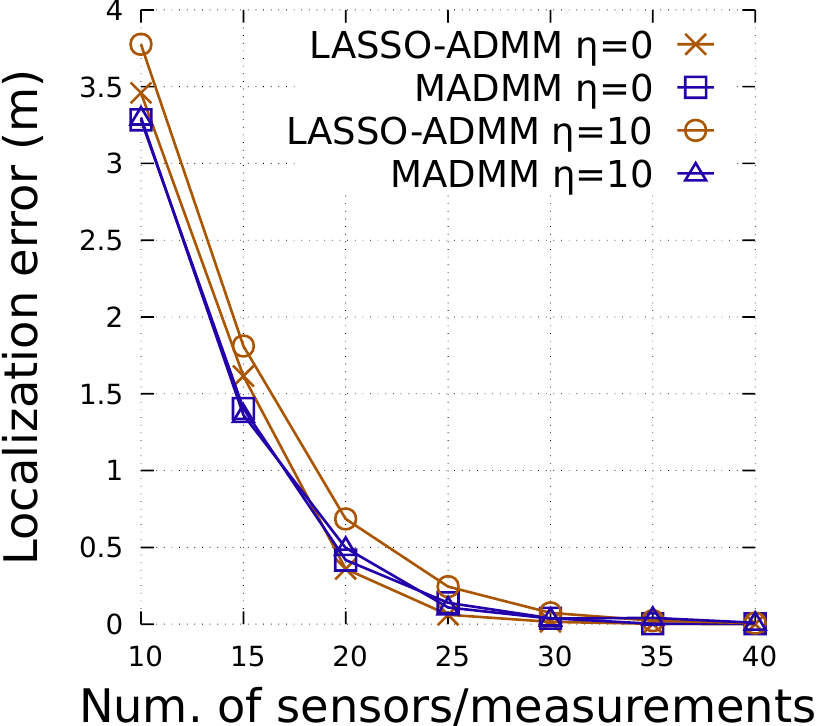}
\includegraphics[width=0.24\textwidth]{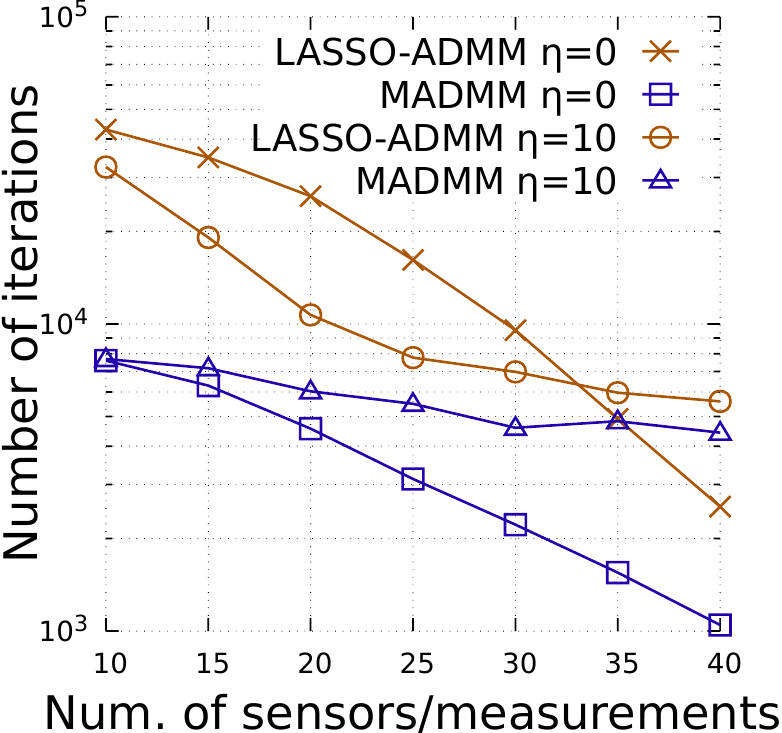}
\caption{Localization error and corresponding number of iterations to converge, varying $m$, with no additive measurement noise and with Gaussian Noise $~\mathcal{N}(0,\eta^2)$, $\eta=10$.}\label{fig:loc}
\end{figure}
We consider $k=4$ targets. Assuming to know $k$, at the end of the procedure we select the $k$ largest values to estimate the occupied cells. We then compute the mean localization error  $\min\frac{1}{k}\sum_{i=1}^k \|\tau_i-\widehat{\tau}_i \|_2$, where $\tau_i$ and $\widehat{\tau}_i$ respectively are the real and estimated positions of the targets (which are assumed to be in the center of the cells). In Figure \ref{fig:loc}, we show the localization error and the number of iterations to converge. We can appreciate that MADMM gives a slightly lower error in a significantly smaller number of iterations, in particular for lower $m$'s.

\section{Conclusions}
In this paper, we have used concave penalization techniques to recover finite-valued sparse signals, with particular focus on the CS framework. We have theoretically proven that the desired signal is the global minimum of a suitable cost functional in the noise-free case. The same functional has been shown to be robust to (signal and measurement) noise. We have then derived recovery algorithms based on ADMM, whose convergence has been discussed. A method to check if a solution is exact also been shown. Numerical experiments show the efficiency of the proposed method with respect to the state-of-the-art in terms of accuracy and convergence speed.

\bibliographystyle{plain}
\bibliography{refs}

\begin{thebibliography}{10}
\providecommand{\url}[1]{#1}
\csname url@samestyle\endcsname
\providecommand{\newblock}{\relax}
\providecommand{\bibinfo}[2]{#2}
\providecommand{\BIBentrySTDinterwordspacing}{\spaceskip=0pt\relax}
\providecommand{\BIBentryALTinterwordstretchfactor}{4}
\providecommand{\BIBentryALTinterwordspacing}{\spaceskip=\fontdimen2\font plus
\BIBentryALTinterwordstretchfactor\fontdimen3\font minus
  \fontdimen4\font\relax}
\providecommand{\BIBforeignlanguage}[2]{{%
\expandafter\ifx\csname l@#1\endcsname\relax
\typeout{** WARNING: IEEEtran.bst: No hyphenation pattern has been}%
\typeout{** loaded for the language `#1'. Using the pattern for}%
\typeout{** the default language instead.}%
\else
\language=\csname l@#1\endcsname
\fi
#2}}
\providecommand{\BIBdecl}{\relax}
\BIBdecl

\bibitem{don06}
D.~L. Donoho, ``Compressed sensing,'' \emph{IEEE Trans. Inf. Theory}, vol.~52,
  no.~4, pp. 1289--1306, 2006.

\bibitem{fou11}
S.~Foucart, ``Hard thresholding pursuit: An algorithm for compressive
  sensing,'' \emph{SIAM J. Numer. Anal.}, vol.~49, no.~6, pp. 2543--2563, 2011.

\bibitem{bio14}
V.~Bioglio, G.~Coluccia, and E.~Magli, ``Sparse image recovery using compressed
  sensing over finite alphabets,'' in \emph{IEEE Int. Conf. Image Process.
  (ICIP)}, 2014, pp. 1287--1291.

\bibitem{bio14sec}
V.~Bioglio, T.~Bianchi, and E.~Magli, ``Secure compressed sensing over finite
  fields,'' in \emph{IEEE Int. Work. Inf. Forensics Secur. (WIFS)}, 2014, pp.
  191--196.

\bibitem{spa15}
S.~Sparrer and R.~F.~H. Fischer, ``Soft-feedback {OMP} for the recovery of
  discrete-valued sparse signals,'' in \emph{European Signal Process. Conf.
  (EUSIPCO)}, 2015, pp. 1461--1465.

\bibitem{ili12}
J.~Ilic and T.~Strohmer, ``Sparsity enhanced decision feedback equalization,''
  \emph{IEEE Trans. Signal Process.}, vol.~60, no.~5, pp. 2422--2432, 2012.

\bibitem{bem99}
A.~Bemporad and M.~Morari, ``Control of systems integrating logic, dynamics,
  and constraints,'' \emph{Automatica}, vol.~35, no.~3, pp. 407--427, 1999.

\bibitem{fen09}
C.~Feng, S.~Valaee, and Z.~Tan, ``Multiple target localization using
  compressive sensing,'' in \emph{IEEE Global Telec. Conf. (GLOBECOM)}, 2009,
  pp. 1--6.

\bibitem{bay15}
A.~Bay, D.~Carrera, S.~M. Fosson, P.~Fragneto, M.~Grella, C.~Ravazzi, and
  E.~Magli, ``Block-sparsity-based localization in wireless sensor networks,''
  \emph{EURASIP J. Wirel. Commun. Netw.}, vol. 2015, no. 182, pp. 1--15, 2015.

\bibitem{baz10}
J.~A. Bazerque and G.~B. Giannakis, ``Distributed spectrum sensing for
  cognitive radio networks by exploiting sparsity,'' \emph{IEEE Trans. Signal
  Process.}, vol.~58, no.~3, pp. 1847--1862, 2010.

\bibitem{zen11}
F.~Zeng, C.~Li, and Z.~Tian, ``Distributed compressive spectrum sensing in
  cooperative multihop cognitive networks,'' \emph{IEEE J. Sel. Top. Sign.
  Process.}, vol.~5, no.~1, pp. 37--48, 2011.

\bibitem{axe12}
E.~Axell, G.~Leus, E.~G. Larsson, and H.~V. Poor, ``Spectrum sensing for
  cognitive radio : State-of-the-art and recent advances,'' \emph{IEEE Signal
  Process. Mag.}, vol.~29, no.~3, pp. 101--116, 2012.

\bibitem{rom13}
D.~Romero and G.~Leus, ``Wideband spectrum sensing from compressed measurements
  using spectral prior information,'' \emph{IEEE Trans. Signal Process.},
  vol.~61, no.~24, pp. 6232--6246, 2013.

\bibitem{das13}
A.~K. Das and S.~Vishwanath, ``On finite alphabet compressive sensing,'' in
  \emph{Proc. IEEE Int. Conf. Acoust. Speech Signal Process. (ICASSP)}, 2013,
  pp. 5890--5894.

\bibitem{fli18}
A.~Flinth and G.~Kutyniok, ``{PROMP}: A sparse recovery approach to
  lattice-valued signals,'' \emph{Appl. Comput. Harmon. Anal.}, vol.~45, no.~3,
  pp. 668--708, 2018.

\bibitem{tia09}
Z.~Tian, G.~Leus, and V.~Lottici, ``Detection of sparse signals under
  finite-alphabet constraints,'' in \emph{IEEE Int. Conf. Acoust. Speech Signal
  Process. (ICASSP)}, 2009, pp. 2349--2352.

\bibitem{sto10}
M.~Stojnic, ``Recovery thresholds for $\ell_1$ optimization in binary
  compressed sensing,'' in \emph{Proc. IEEE Int. Symp. Inf. Theory (ISIT)},
  2010, pp. 1593--1597.

\bibitem{lee16}
N.~Lee, ``{MAP} support detection for greedy sparse signal recovery algorithms
  in compressive sensing,'' \emph{IEEE Trans. Signal Process.}, vol.~64,
  no.~19, pp. 4987--4999, 2016.

\bibitem{fox18asi}
S.~M. Fosson, ``Non-convex approach to binary compressed sensing,'' in
  \emph{Asilomar Conf. Signals Syst. Comput.}, 2018.

\bibitem{kei17}
S.~Keiper, G.~Kutyniok, D.~G. Lee, and G.~E. Pfander, ``Compressed sensing for
  finite-valued signals,'' \emph{Linear Algebra and its Applications}, vol.
  532, no. Supplement C, pp. 570--613, 2017.

\bibitem{can08rew}
E.~J. Cand\`es, M.~B. Wakin, and S.~Boyd, ``Enhancing sparsity by reweighted
  $\ell_1$ minimization,'' \emph{J. Fourier Anal. Appl.}, vol.~14, no. 5-6, pp.
  877--905, 2008.

\bibitem{gas09}
G.~Gasso, A.~Rakotomamonjy, and S.~Canu, ``Recovering sparse signals with a
  certain family of nonconvex penalties and dc programming,'' \emph{IEEE Trans.
  Signal Process.}, vol.~57, no.~12, pp. 4686--4698, 2009.

\bibitem{woo16}
J.~Woodworth and R.~Chartrand, ``Compressed sensing recovery via nonconvex
  shrinkage penalties,'' \emph{Inverse Problems}, vol.~32, no.~7, pp.
  75\,004--75\,028, 2016.

\bibitem{cha13}
R.~Chartrand, E.~Y. Sidky, and X.~Pan, ``Nonconvex compressive sensing for
  x-ray ct: An algorithm comparison,'' in \emph{Asilomar Conf. Signals Syst.
  Comput.}, 2013, pp. 665--669.

\bibitem{cha09}
R.~Chartrand, ``Fast algorithms for nonconvex compressive sensing: Mri
  reconstruction from very few data,'' in \emph{IEEE Int. Symp. Biomed. Imag.},
  2009, pp. 262--265.

\bibitem{lap12}
L.~Laporte, R.~Flamary, S.~Canu, S.~D\'{e}jean, and J.~Mothe, ``Nonconvex
  regularizations for feature selection in ranking with sparse {SVM},''
  \emph{IEEE Trans. Neural Netw. Learn. Syst.}, vol.~25, no.~6, pp. 1118--1130,
  2014.

\bibitem{tib96}
R.~Tibshirani, ``Regression shrinkage and selection via the lasso,''
  \emph{Journal of the Royal Statistical Society, Series B}, vol.~58, pp.
  267--288, 1996.

\bibitem{dau04}
I.~Daubechies, M.~Defrise, and C.~De~Mol, ``An iterative thresholding algorithm
  for linear inverse problems with a sparsity constraint,'' \emph{Commun. Pure
  Appl. Math.}, vol.~57, no.~11, pp. 1413 -- 1457, 2004.

\bibitem{for10}
M.~Fornasier, ``Numerical methods for sparse recovery,'' in \emph{Theoretical
  Foundations and Numerical Methods for Sparse Recovery}, M.~Fornasier,
  Ed.\hskip 1em plus 0.5em minus 0.4em\relax Radon Series Comp. Appl. Math., de
  Gruyter, 2010, pp. 93--200.

\bibitem{boy10}
S.~Boyd, N.~Parikh, E.~Chu, B.~Peleato, and J.~Eckstein, ``Distributed
  optimization and statistical learning via the alternating direction method of
  multipliers,'' \emph{Found. Trends Mach. Learn.}, vol.~3, no.~1, pp. 1 --
  122, 2010.

\bibitem{mata15}
J.~Matamoros, S.~M. Fosson, E.~Magli, and C.~Ant\'{o}n-Haro, ``Distributed
  {ADMM} for in-network reconstruction of sparse signals with innovations,''
  \emph{IEEE Trans. Signal Inf. Process. Netw.}, vol.~1, no.~4, pp. 225--234,
  2015.

\bibitem{fia18}
A.~Fiandrotti, S.~M. Fosson, C.~Ravazzi, and E.~Magli, ``Gpu-accelerated
  algorithms for compressed signals recovery with application to astronomical
  imagery deblurring,'' \emph{Int. J. Remote Sens.}, vol.~39, no.~7, pp.
  2043--2065, 2018.

\bibitem{zha10MCP}
C.-H. Zhang, ``Nearly unbiased variable selection under minimax concave
  penalty,'' \emph{Ann. Statist.}, vol.~38, no.~2, pp. 894--942, 2010.

\bibitem{fan01_pioneer}
J.~Fan and R.~Li, ``Variable selection via nonconcave penalized likelihood and
  its oracle properties,'' \emph{J. Amer. Statist. Assoc.}, vol.~96, no. 456,
  pp. 1348--1360, 2001.

\bibitem{fan11}
J.~Fan and J.~Lv, ``Nonconcave penalized likelihood with np-dimensionality,''
  \emph{IEEE Trans. Inf. Theory}, vol.~57, no.~8, pp. 5467--5484, 2011.

\bibitem{fan14}
J.~Fan, L.~Xue, and H.~Zou, ``Strong oracle optimality of folded concave
  penalized estimation,'' \emph{Ann. Statist.}, vol.~42, no.~3, pp. 819--849,
  2014.

\bibitem{zou08_LLA}
H.~Zou and R.~Li, ``One-step sparse estimates in nonconcave penalized
  likelihood models,'' \emph{Annals of Statistics}, vol.~36, no.~4, p. 1509,
  2008.

\bibitem{zha12}
C.-H. Zhang and T.~Zhang, ``A general theory of concave regularization for
  high-dimensional sparse estimation problems,'' \emph{Statist. Sci.}, vol.~27,
  no.~4, pp. 576 --593, 2012.

\bibitem{irls}
I.~Daubechies, R.~DeVore, M.~Fornasier, and C.~S. Gunturk, ``Iteratively
  reweighted least squares minimization for sparse recovery,'' \emph{Commun.
  Pure Appli. Math.}, vol.~63, no.~1, pp. 1--38, 2010.

\bibitem{rav15irls}
C.~Ravazzi and E.~Magli, ``Gaussian mixtures based {IRLS} for sparse recovery
  with quadratic convergence,'' \emph{IEEE Trans. Signal Process.}, vol.~63,
  no.~13, pp. 3474--3489, 2015.

\bibitem{faz03}
M.~Fazel, H.~Hindi, and S.~Boyd, ``Log-det heuristic for matrix rank
  minimization with applications to {H}ankel and {E}uclidean distance
  matrices,'' in \emph{IEEE Proc. American Control Conf. (ACC)}, vol.~3, 2003,
  pp. 2156--2162.

\bibitem{cal16}
M.~Calvo-Fullana, J.~Matamoros, C.~Ant\'{o}n-Haro, and S.~M. Fosson,
  ``Sparsity-promoting sensor selection with energy harvesting constraints,''
  in \emph{Proc. IEEE Int. Conf. Acoust. Speech Signal Process. (ICASSP)},
  2016, pp. 3766--3770.

\bibitem{fox16}
S.~M. Fosson, J.~Matamoros, C.~Ant\'{o}n-Haro, and E.~Magli, ``Distributed
  recovery of jointly sparse signals under communication constraints,''
  \emph{IEEE Trans. Signal Process.}, vol.~64, no.~13, pp. 3470--3482, 2016.

\bibitem{hua18}
X.~Huang and M.~Yan, ``Nonconvex penalties with analytical solutions for
  one-bit compressive sensing,'' \emph{Signal Process.}, vol. 144, pp.
  341--351, 2018.

\bibitem{can08}
E.~J. Cand\`es, ``The restricted isometry property and its implications for
  compressed sensing,'' \emph{C. R. Acad. Sci. Paris, Ser. I}, vol. 346, pp.
  589--592, 2008.

\bibitem{las15}
J.~B. Lasserre, \emph{An Introduction to Polynomial and Semi-Algebraic
  Optimization}, ser. Cambridge Texts in Applied Mathematics.\hskip 1em plus
  0.5em minus 0.4em\relax Cambridge University Press, 2015.

\bibitem{fou13}
S.~Foucart and H.~Rauhut, \emph{A Mathematical Introduction to Compressive
  Sensing}.\hskip 1em plus 0.5em minus 0.4em\relax New York: Springer, 2013.

\bibitem{matrix}
R.~A. Horn and C.~R. Johnson, \emph{Matrix Analysis - II Ed.}\hskip 1em plus
  0.5em minus 0.4em\relax Cambridge University Press, 2013.

\bibitem{hon15}
M.~Hong, Z.~Q. Luo, and M.~Razaviyayn, ``Convergence analysis of alternating
  direction method of multipliers for a family of nonconvex problems,'' in
  \emph{2015 IEEE Int. Conf. Acoust. Speech Signal Process. (ICASSP)}, 2015,
  pp. 3836--3840.

\bibitem{wan15}
Y.~Wang, W.~Yin, and J.~Zeng, ``Global convergence of {ADMM} in nonconvex
  nonsmooth optimization,'' \emph{J. Scientific Comput.}, 2018.

\bibitem{li15}
G.~Li and T.~K. Pong, ``Global convergence of splitting methods for nonconvex
  composite optimization,'' \emph{SIAM J. Optim.}, vol.~25, no.~4, pp.
  2434--2460, 2015.

\bibitem{hon16}
M.~Hong, Z.~Q. Luo, and M.~Razaviyayn, ``Convergence analysis of alternating
  direction method of multipliers for a family of nonconvex problems,''
  \emph{SIAM J. Optim.}, vol.~26, no.~1, pp. 337--364, 2016.

\bibitem{tib13}
R.~J. Tibshirani, ``{The {L}asso problem and uniqueness},'' \emph{Electronic
  Journal of Statistics}, vol.~7, pp. 1456--1490, 2013.

\bibitem{rfm15}
C.~Ravazzi, S.~M. Fosson, and E.~Magli, ``Distributed iterative thresholding
  for $\ell_0$/$\ell_1$-regularized linear inverse problems,'' \emph{IEEE
  Trans. Inf. Theory}, vol.~61, no.~4, pp. 2081--2100, 2015.

\bibitem{bayram16}
Ä.~Bayram, ``On the convergence of the iterative shrinkage/thresholding
  algorithm with a weakly convex penalty,'' \emph{IEEE Trans. Signal Process.},
  vol.~64, no.~6, pp. 1597--1608, 2016.

\end{thebibliography}

\end{document}